\documentclass[11pt]{article}

\usepackage{amssymb, amsfonts, amsmath, amsthm, graphics}
\usepackage{lineno}

\newtheorem{theorem}{Theorem}[section]
\newtheorem{lemma}[theorem]{Lemma}
\newtheorem{proposition}[theorem]{Proposition}
\newtheorem{corollary}[theorem]{Corollary}

\theoremstyle{definition}
\newtheorem{definition}[theorem]{Definition}
\newtheorem{notation}[theorem]{Notation}
\newtheorem{remark}[theorem]{Remark}
\newtheorem{example}[theorem]{Example}

\newcommand{\cA}{ {\cal A} }
\newcommand{\cB}{ {\cal B} }
\newcommand{\bC}{ {\mathbb C} }

\newcommand{\cD}{ {\cal D} }
\newcommand{\uE}{ {\mathbf e} }
\newcommand{\cG}{ {\cal G} }
\newcommand{\uH}{ {\mathbf h} }
\newcommand{\bI}{ {\mathbb I} }
\newcommand{\cM}{{\cal M}}

\newcommand{\uP}{ {\mathbf p} }
\newcommand{\cR}{{\cal R}}

\newcommand{\uU}{ {\mathbf u} }
\newcommand{\bX}{ {\mathbb X} }

\newcommand{\cYk}{ {\cal Y}^{(k)} }
\newcommand{\cYone}{ {\cal Y}^{(1)} }

\newcommand{\dalg}{ \cD_{\mathrm{alg}} }
\newcommand{\cf}{ \mbox{Cf} }
\newcommand{\cigma}{ \Gamma }
\newcommand{\id}{ \mbox{id} }
\newcommand{\ls}{ LS }
\newcommand{\mult}{ \mbox{Mult} }
\newcommand{\sym}{ \mbox{Sym} }
\newcommand{\term}{ \mbox{term} }
\newcommand{\PrimY}{ \mbox{Prim} \bigl( \, \cYone \, \bigr) }

\newcommand{\ee}{\varepsilon}

\newcommand{\ens}{ ( e_n )_{n=1}^{\infty} }
\newcommand{\hns}{ ( h_n )_{n=1}^{\infty} }
\newcommand{\pns}{ ( p_n )_{n=1}^{\infty} }
\newcommand{\sns}{ ( s_n )_{n=1}^{\infty} }
\newcommand{\tns}{ ( t_n )_{n=1}^{\infty} }

\begin{document}

$\ $

\begin{center}
{\bf\Large Hopf algebras and the logarithm of the S-transform}

\vspace{6pt}

{\bf\Large in free probability}

\vspace{20pt}

{\large Mitja Mastnak \hspace{2cm}
Alexandru Nica \footnote{Research supported by a Discovery Grant
from NSERC, Canada.} }

\vspace{10pt}

\end{center}

\begin{abstract}

\noindent
Let $k$ be a positive integer and let $\cG_k$ denote the set
of all joint distributions of $k$-tuples $(a_1, \ldots , a_k)$
in a non-commutative probability space $( \cA , \varphi )$
such that $\varphi (a_1) = \cdots = \varphi (a_k) = 1$.
$\cG_k$ is a group under the operation of free multiplicative
convolution $\boxtimes$. We identify
$\bigl( \, \cG_k, \boxtimes \, \bigr)$ as the group of characters
of a certain Hopf algebra $\cYk$. Then, by using the log map
from characters to infinitesimal characters of $\cYk$, we
introduce a transform $\ls_{\mu}$ for distributions $\mu \in \cG_k$.
$\ls_{\mu}$ is a power series in $k$ non-commuting indeterminates
$z_1, \ldots , z_k$; its coefficients can be computed from the
coefficients of the $R$-transform of $\mu$ by using summations
over chains in the lattices $NC(n)$ of non-crossing partitions.
The $LS$-transform has the ``linearizing'' property that
\[
\ls_{\mu \boxtimes \nu} = \ls_{\mu} + \ls_{\nu}, \ \ \forall \,
\mu , \nu \in \cG_k \mbox{ such that }
\mu \boxtimes \nu = \nu \boxtimes \mu .
\]

In the particular case $k=1$ one has that $\cYone$ is naturally
isomorphic to the Hopf algebra $\sym$ of symmetric functions, and
that the $LS$-transform is very closely related to the logarithm
of the $S$-transform of Voiculescu, by the formula
\[
\ls_{\mu} (z) = -z \log S_{\mu} (z), \ \ \forall \,
\mu \in \cG_1.
\]
In this case the group $( \cG_1, \boxtimes )$ can be identified
as the group of characters of $\sym$, in such a way that the
$S$-transform, its reciprocal $1/S$ and its logarithm
$\log S$ relate in a natural sense to the sequences of complete,
elementary and respectively power sum symmetric functions.
\end{abstract}

\vspace{6pt}

\begin{center}
{\bf\large 1. Introduction}
\end{center}
\setcounter{section}{1}

\noindent
In this paper we study joint distributions of $k$-tuples of elements
in a noncommutative probability space. Let $( \cA , \varphi )$ be
a noncommutative probability space (i.e. $\cA$ is a unital algebra
over $\bC$ and $\varphi : \cA \to \bC$ is a linear functional such
that $\varphi ( 1_{\cA} ) = 1$), and let $a_1, \ldots, a_k$ be
elements of $\cA$. The distribution of $(a_1, \ldots , a_k)$
is the linear functional $\mu$ on the algebra of noncommutative
polynomials $\bC \langle X_1, \ldots , X_k \rangle$ defined by the
requirement that 
\[
\mu ( X_{i_1} \cdots X_{i_n} ) = \varphi ( a_{i_1} \cdots a_{i_n} ),
\ \ \forall \, n \geq 0, \ \forall \, 1 \leq i_1, \ldots , i_n \leq k.
\]
We denote by $\dalg (k)$ the set of linear functionals on
$\bC \langle X_1, \ldots , X_k \rangle$ that arise in this way.
(Clearly, this is just the set of all linear functionals on
$\bC \langle X_1, \ldots , X_k \rangle$ such that $\mu (1) = 1$.)
On $\dalg (k)$ we have a binary operation
$\boxtimes$ which reflects the multiplication of two freely
independent $k$-tuples in a noncommutative probability space.
That is, $\boxtimes$ is well-defined and uniquely determined by the
following requirement: if $a_1, \ldots , a_k, b_1, \ldots , b_k$ are
elements in a noncommutative probability space $( \cA , \varphi )$
such that $(a_1, \ldots , a_k)$ has distribution $\mu$,
$(b_1, \ldots , b_k)$ has distribution $\nu$,
and $\{ a_1, \ldots , a_k \}$ is freely independent from
$\{ b_1, \ldots , b_k \}$, then it follows that the distribution
of $( a_1 b_1, \ldots , a_k b_k )$ is equal to $\mu \boxtimes \nu$.
The operation $\boxtimes$ on $\dalg (k)$ is associative and unital,
where the unit is the functional $\mu_o \in \dalg (k)$ with
$\mu_o ( X_{i_1} \cdots X_{i_n} )= 1$ for all $n \geq 1$ and
$1 \leq i_1, \ldots , i_n \leq k$. A distribution $\mu \in \dalg (k)$
is invertible with respect to $\boxtimes$ if and only if it satisfies
$\mu (X_i) \neq 0$, $\forall \, 1 \leq i \leq k$; and moreover, the
subset
\begin{equation}  \label{eqn:1.1}
\cG_k := \{ \mu \in \dalg (k) \mid \mu (X_i)= 1, \ \
\forall \, 1 \leq i \leq k \}
\end{equation}
is a subgroup in the group of invertibles with respect to $\boxtimes$.
For a basic introduction to free multiplicative convolution, we refer
to Section 3.6 of \cite{VDN92} or to Lecture 14 in \cite{NS06}.

The main goal of the present paper is to introduce a transform
$\ls_{\mu}$ for distributions $\mu \in \cG_k$, which linearizes
commuting $\boxtimes$-products. We arrive to the $LS$-transform
by identifying $\bigl( \, \cG_k , \boxtimes \, \bigr)$ as the
group of characters of a certain Hopf algebra $\cYk$. As an
algebra, $\cYk$ is merely a commutative algebra of polynomials:
\begin{equation}  \label{eqn:1.2}
\cYk = \bC \bigl[ \, Y_w \mid w \in [k]^*, \ |w| \geq 2 \, \bigr],
\end{equation}
where $[k]^*$ is the set of all words of finite length made with
letters from the alphabet $\{ 1, \ldots , k \}$, and where
$|w|$ denotes the length (i.e. number of letters) of a word
$w \in [k]^*$. The feature which relates $\cYk$ to free probability
is its comultiplication: for $w \in [k]^*$ with $|w| =n \geq 2$, the
comultiplication $\Delta (Y_w) \in \cYk \otimes \cYk$ is defined by
using a special type of summation over the lattice of non-crossing
partitions $NC(n)$ which is known (see e.g. Lectures 14 and 17 of
\cite{NS06}) to relate to the multiplication of free random variables.
The details of the construction of $\cYk$ are presented in
Section 3 below. Here we only mention the fact that $\cYk$ is a
graded connected Hopf algebra, where every generator $Y_w$ from
(\ref{eqn:1.2}) is homogeneous of degree $|w|-1$. A consequence of
this fact which is important for the present paper
is that one can then define the exponential $\exp \xi$ for every
linear functional $\xi: \cYk \to \bC$ such that
$\xi ( 1 ) = 0$, and one can define the logarithm $\log \eta$
for every linear functional $\eta : \cYk \to \bC$ such that
$\eta ( 1 ) = 1$. Moreover, the maps $\exp$ and $\log$ 
are bijections inverse to each other, between
the two sets of linear functionals mentioned above.

Let $\bX ( \cYk )$ be the set of characters of $\cYk$,
\begin{equation}  \label{eqn:1.3}
\bX ( \cYk ) := \Bigl\{ \eta : \cYk \to \bC \mid
\eta ( 1 ) = 1, \mbox{ $\eta$ is linear and multiplicative}
\Bigr\}.
\end{equation}
Then the logarithm map sends $\bX ( \cYk )$ bijectively onto the
set $\bI ( \cYk )$ of infinitesimal characters of $\cYk$; the
latter set is defined by
\begin{equation}  \label{eqn:1.4}
\bI ( \cYk ) := \left\{ \xi : \cYk \to \bC  \begin{array}{cl}
\vline & \mbox{ $\xi$ is linear and satisfies}  \\
\vline & \xi (PQ) = \xi (P) \ee (Q) + \ee (P) \xi (Q), \ \
         \forall \, P,Q \in \cYk
\end{array}  \right\} ,
\end{equation}
where $\ee$ is the counit of $\cYk$ ($\ee \in \bX ( \cYk )$, and
is uniquely determined by the requirement that
$\ee ( Y_w ) = 0$ for every $w \in [k]^*$ with $|w| \geq 2$).
Let us also record the facts that $\bX ( \cYk )$ is a group under
the operation of convolution for linear 
functionals on $\cYk$, and
that the log map linearizes commuting convolution products:
\begin{equation}  \label{eqn:1.5}
\Bigl( \, \eta_1, \eta_2 \in \bX ( \cYk ),
\ \eta_1  \eta_2 = \eta_2  \eta_1 \, \Bigr)
\Rightarrow  \log ( \eta_1  \eta_2) = \log \eta_1 + \log \eta_2 .
\end{equation}
(In this paper the convolution of linear functionals on a Hopf
algebra is denoted as a plain multiplication. Our conventions of
notation and a few background facts about graded connected Hopf
algebras are collected in Section 2C of the paper.)

So now let us state precisely what is the connection between the
Hopf algebra $\cYk$ and the operation $\boxtimes$. Let $\mu$ be
a distribution in $\cG_k$, and let us consider its $R$-transform
$R_{\mu}$; this is a power series in $k$ non-commuting indeterminates
$z_1, \ldots , z_k$, of the form
\begin{equation}  \label{eqn:1.6}
R_{\mu} ( z_1, \ldots , z_k ) = \sum_{i=1}^k z_i +
\sum_{ \begin{array}{c}
{\scriptstyle w \in [k]^* , } \\
{\scriptstyle |w| \geq 2}
\end{array} } \ \alpha_{w} z_w,
\end{equation}
where in the latter sum we used the shorthand notation
$z_w := z_{i_1} \cdots z_{i_n}$ for
$w = ( i_1, \ldots , i_n ) \in [k]^*$ with $n \geq 2$.
A brief review of how the coefficients $\alpha_w$ are
calculated from the moments (i.e.  values on monomials) of $\mu$
is made in Section 2B below -- see Definition \ref{def:2.5},
Remark \ref{rem:2.6}.

\begin{definition}  \label{def:1.1}
Let $\mu$ be a distribution in $\cG_k$ and consider the
$R$-transform $R_{\mu}$, denoted as in Equation (\ref{eqn:1.6}).
The {\em character of $\cYk$ associated to $\mu$} is the character
$\chi_{\mu} \in \bX ( \cYk )$ uniquely determined by the requirement
that
\begin{equation}  \label{eqn:1.7}
\chi_{\mu} (Y_w) = \alpha_w, \ \ \forall \, w \in [k]^*
\mbox{ such that } |w| \geq 2.
\end{equation}
\end{definition}

\begin{theorem}  \label{thm:1.2}
The map $\mu \mapsto \chi_{\mu}$ defined above is a group
isomorphism from $\bigl( \, \cG_k , \boxtimes \, \bigr)$ onto
the group $\bX ( \cYk )$ of characters of $\cYk$ (endowed
with the operation of convolution).
\end{theorem}

The proof of Theorem \ref{thm:1.2} is presented in Section 3 below.
Let us mention that the proof is quite short, and follows easily
from known facts about how the multiplication of free $k$-tuples
is described in terms of their $R$-transforms, as explained for
instance in Lecture 17 of \cite{NS06}. (In a certain sense, we are
dealing here with a situation where the convolution of characters
of $\cYk$ was studied before looking at $\cYk$ itself.)
Nevertheless, Theorem \ref{thm:1.2} is important because it brings
to attention the fact that Hopf algebra methods can be used in the
study of $\boxtimes$. In particular, by taking (\ref{eqn:1.5}) into
account, we see that Theorem \ref{thm:1.2} has the following
immediate corollary.

\begin{corollary}  \label{cor:1.3}
Let $\mu$ and $\nu$ be distributions in $\cG_k$ such that
$\mu \boxtimes \nu = \nu \boxtimes \mu$. Then
\begin{equation}  \label{eqn:1.8}
\log \chi_{\mu \boxtimes \nu}  =
\log \chi_{\mu} + \log \chi_{\nu} ,
\end{equation}
where the characters $\chi_{\mu} , \ \chi_{\nu}$ and
$\chi_{\mu \boxtimes \nu}$ are as in Definition \ref{def:1.1},
and their logarithms are the corresponding functionals from
$\bI ( \cYk )$.
\end{corollary}

The $LS$-transform $\ls_{\mu}$ is defined so that it stores the
information about the infinitesimal character $\log \chi_{\mu}$,
as follows.

\begin{definition}  \label{def:1.4}
Let $\mu$ be a distribution in $\cG_k$. The $LS$-transform of $\mu$
is the power series
\begin{equation}  \label{eqn:1.9}
\ls_{\mu} (z_1, \ldots , z_k) :=
\sum_{ \begin{array}{c}
{\scriptstyle  w \in [k]^* , }  \\
{\scriptstyle |w| \geq 2  }
\end{array} } \
\Bigl( \, \bigl( \log \chi_{\mu} \bigr) (Y_w) \, \Bigr) z_w,
\end{equation}
where $\log \chi_{\mu} \in \bI ( \cYk )$ is as in
Corollary \ref{cor:1.3}, and where the meaning of ``$z_w$'' is
same as in Equation (\ref{eqn:1.6}).
\end{definition}

Clearly, Corollary \ref{cor:1.3} can also be phrased as a statement
about $LS$-transforms. In this guise, it simply says that the
$LS$-transform linearizes commuting $\boxtimes$-products:

\begin{corollary}  \label{cor:1.5}
Let $\mu$ and $\nu$ be distributions in $\cG_k$ such that
$\mu \boxtimes \nu = \nu \boxtimes \mu$. Then
\begin{equation}  \label{eqn:1.95}
\ls_{\mu \boxtimes \nu}  = \ls_{\mu} + \ls_{\nu} .
\end{equation}
\end{corollary}

In particular, formula (\ref{eqn:1.95}) always applies when one
of $\mu$, $\nu$ is the joint
distribution of a repeated $k$-tuple $(a,a, \ldots , a)$, where
$a$ is a random variable in a non-commutative probability space
$( \cA , \varphi )$. A version of this fact which lives in the 
framework of a $C^*$-probability space is discussed in Example 
\ref{ex:5.2} below.

The $LS$-transform was introduced above by using the
Hopf algebra $\cYk$, but it can also be described directly in
combinatorial terms, by using summations over chains in lattices
of non-crossing partitions. We next explain how this goes.

A chain in $NC(n)$ is an object of the form
\begin{equation}  \label{eqn:1.10}
\cigma = ( \pi_0, \pi_1, \ldots , \pi_{\ell} )
\end{equation}
with $\pi_0, \pi_1, \ldots , \pi_{\ell} \in NC(n)$ such that
$0_n = \pi_0< \pi_1< \cdots < \pi_{\ell} = 1_n$ (and where $0_n$
and $1_n$ denote the minimal and maximal element of $NC(n)$,
respectively). For a chain $\cigma$ as in (\ref{eqn:1.10}), the
number $\ell$ is called the length of $\cigma$ and will be
denoted as $| \cigma |$. Given a formal power series
$f$ in non-commuting indeterminates $z_1, \ldots , z_k$,
one has a natural way of defining some ``generalized coefficients''
\begin{equation}  \label{eqn:1.11}
\cf_{(i_1, \ldots , i_n)}^{ \, ( \cigma )} (f)
\end{equation}
where $n \geq 2$, $1 \leq i_1, \ldots , i_n \leq k$, and
$\cigma$ is a chain in $NC(n)$. Every generalized
coefficient (\ref{eqn:1.11}) is defined to be a certain product of
actual coefficients of $f$ (see Definition \ref{def:4.3} below for
the precise formula).
By using the generalized coefficients (\ref{eqn:1.11}), the
combinatorial description of the $LS$-transform is stated as follows.

\begin{theorem}  \label{thm:1.6}
Let $\mu$ be a distribution in $\cG_k$, and let
$w = ( i_1, \ldots , i_n)$ be a word in $[k]^*$, where $n \geq 2$.
Consider (as in Equation (\ref{eqn:1.9}) of Definition \ref{def:1.4})
the coefficient $\bigl( \log \chi_{\mu} \bigr) (Y_w)$ of $z_w$ in the
$LS$-transform of $\mu$. This coefficient can be also expressed as
\begin{equation}  \label{eqn:1.12}
\bigl( \log \chi_{\mu} \bigr) (Y_w) = \sum_{ \begin{array}{c}
{\scriptstyle \cigma \ \mathrm{chain} } \\
{\scriptstyle \mathrm{in} \ NC(n) }
\end{array}  }  \ \frac{ (-1)^{1+| \cigma |} }{| \cigma |}
\cf_{(i_1, \ldots , i_n)}^{ \, ( \cigma )} ( R_{\mu} ),
\end{equation}
where $R_{\mu}$ is the $R$-transform of $\mu$.
\end{theorem}

\begin{remark}  \label{rem:1.7}

$1^o$ $NC(n)$ has a unique chain of length 1 (namely the chain
$(0_n, 1_n)$); the term indexed by this chain in the sum on the
right-hand side of (\ref{eqn:1.12}) is precisely the coefficient
of $z_w$ in $R_{\mu}$, while every other term of the same sum
turns out to be a product of coefficients of $R_{\mu}$ with
lengths strictly smaller than $n$. From this observation it is
immediate that the coefficients of $R_{\mu}$ can be computed back,
recursively, in terms of the coefficients of $\ls_{\mu}$. Moreover,
since $\mu$ is completely determined by its $R$-transform, we thus
see that $\mu$ is completely determined by the series $\ls_{\mu}$
as well.

$2^o$ By using Theorem \ref{thm:1.6} one finds that the
$LS$-transform inherits a fundamental property that the 
$R$-transform is known to have in connection to free independence:
a $k$-tuple in a noncommutative probability
space is freely independent if and only if the $LS$-transform of
its distribution separates the variables. The precise statement of
this fact appears as Proposition \ref{prop:5.4} below.
\end{remark}

\vspace{10pt}

In the remaining part of the introduction we look at the particular
case when $k=1$. In this case the notations are simplified due to the
fact that words over the 1-letter alphabet $\{ 1 \}$ are determined
by their lengths. We make the convention to write simply ``$Y_n$''
instead of $Y_{(1,1, \ldots ,1)}$ with $n$ repetitions of 1 in the
index; thus Equation (\ref{eqn:1.2}) is now written in the form
\begin{equation}  \label{eqn:1.13}
\cYone = \bC \bigl[ \, Y_n \mid n \geq 2 \, \bigr],
\end{equation}
while Equation (\ref{eqn:1.9}) defining $\ls_{\mu}$ reduces to
\begin{equation}  \label{eqn:1.14}
\ls_{\mu} (z) = \sum_{n=2}^{\infty} \Bigl( \, ( \log \chi_{\mu} )
(Y_n) \, \Bigr) z^n.
\end{equation}

A special feature of the case $k=1$ (not holding for $k \geq 2$) is 
that the operation $\boxtimes$ is commutative. Hence for $k=1$ the 
linearization property stated in Corollary
\ref{cor:1.5} holds for all $\mu, \nu \in \cG_1$.

Now, in the case $k=1$ there exists an established way of treating
multiplicative free convolution, by using Voiculescu's
{\em $S$-transform}. The $S$-transform $S_{\mu}$ of a distribution
$\mu \in \cG_1$ is a power series in an indeterminate $z$, with
constant term equal to 1. Taking $S$-transforms converts $\boxtimes$
into plain multiplication of power series:
\begin{equation}  \label{eqn:1.15}
S_{\mu \boxtimes \nu} (z) = S_{\mu} (z) \cdot S_{\nu} (z), \ \
\forall \, \mu , \nu \in \cG_1
\end{equation}
(see Section 3.6 in \cite{VDN92}, or Lecture 18 in \cite{NS06}).
We prove that the one-dimensional $\ls$-transform is related to
the $S$-transform, as follows.

\begin{theorem}  \label{thm:1.8}
For a distribution $\mu \in \cG_1$, the power series
$S_{\mu}$ and $\ls_{\mu}$ are related by
\begin{equation}  \label{eqn:1.17}
\ls_{\mu} (z) = -z \log S_{\mu} (z).
\end{equation}
\end{theorem}

\vspace{6pt}

Theorem \ref{thm:1.8} shows that the multiplicativity of the 
$S$-transform can be retrieved from the particular case $k=1$ of
Corollary \ref{cor:1.5}: one only needs to divide by $-z$ and then 
exponentiate both sides of Equation (\ref{eqn:1.95}). We should note 
here that this exponentiating trick is specific to the 1-variable 
situation, for $k \geq 2$ it is probably better to work with 
$\ls_{\mu}$ itself rather than considering its exponential. (Indeed, 
for $k \geq 2$ it isn't generally true that the
$k$-variable series $\ls_{\mu}$ and $\ls_{\nu}$ would commute,
even if $\mu , \nu \in \cG_k$ are such that
$\mu \boxtimes \nu = \nu \boxtimes \mu$. Hence when one
exponentiates Equation (\ref{eqn:1.95}) for $k \geq 2$, the
series appearing on the right-hand side isn't generally equal to
$\exp \bigl( \ls_{\mu} \bigr) \cdot \exp \bigl( \ls_{\nu} \bigr)$.)

The proof of Theorem \ref{thm:1.8} is obtained by following
the connections that $\cYone$ has with symmetric functions. We start
from the fact that (independently of the
isomorphism $\cG_1 \ni \mu \mapsto \chi_{\mu} \in \bX ( \cYone )$
from Theorem \ref{thm:1.2}) one can also find a natural
group isomorphism $\mu \mapsto \theta_{\mu}$ from
$( \cG_1, \boxtimes )$ onto the group of characters of the Hopf 
algebra $\sym$ of symmetric functions.
The map $\mu \mapsto \theta_{\mu}$ is defined in such a way that
the $S$-transform, its reciprocal $1/S$ and its logarithm
$\log S$ relate in a natural sense to the sequences of complete,
elementary and respectively power sum symmetric functions (see
Remark \ref{rem:6.7} below for the precise description of how
this happens). We then concretely put into evidence an isomorphism
$\Phi : \cYone \to \sym$ which has the property that
\begin{equation}   \label{eqn:1.18}
\chi_{\mu} = \theta_{\mu} \circ \Phi, \ \ 
\forall \, \mu \in \cG_1;
\end{equation}
with the help of $\Phi$ we can place the 1-dimensional $LS$-transform
in the framework of $\sym$, and then prove the relation stated in 
Theorem \ref{thm:1.8} (this is done in Proposition \ref{prop:6.10}
and Corollary \ref{cor:6.12}).

It is interesting to note that the way the $\ls$-transform is
introduced in this paper is close in spirit to how the
$S$-transform was first found by Voiculescu in \cite{V87}. Some
other approaches to the $S$-transform were found in the meanwhile
(e.g. in \cite{H97} or in \cite{NS97}), but the one
from \cite{V87} has the distinctive feature that it relies on the
exponential map for a certain commutative Lie group, constructed to
reflect how the moments of $\mu \boxtimes \nu$ are computed in
terms of the moments of $\mu$ and of $\nu$
($\mu , \nu \in \cG_1$). The approach in the present paper goes on
the same lines, with the difference that it insists on structures
(free cumulants, Hopf algebras) where it is easier to pursue a 
detailed combinatorial analysis. The possibility of moving up to
$\ls$-transforms for $k$-tuples then comes as a manifestation of 
the general principle that combinatorial arguments in free 
probability often extend without much trouble from the 1-variable
to the $k$-variable setting.

We conclude this introduction by describing how the paper is
organized. Besides the introduction, the paper has six other sections.
Section 2 contains a review of some background and notations.
In Section 3 we introduce the Hopf algebra $\cYk$ and we prove
Theorem \ref{thm:1.2}, then Section 4 is devoted to the combinatorics
of the $\ls$-transform and to the proof of Theorem \ref{thm:1.6}.
In Section 5 we discuss some basic properties of $\ls_{\mu}$.
The last two sections of the paper are devoted to the one-variable 
framework: in Section 6 we prove Theorem \ref{thm:1.8} and 
in Section 7 we discuss in more detail the isomorphism between
$\cYone$ and the Hopf algebra $\sym$ of symmetric functions. 

$\ $

$\ $

\begin{center}
{\bf\large 2. Background and notations}
\end{center}
\setcounter{section}{2}
\setcounter{equation}{0}
\setcounter{theorem}{0}

\vspace{4pt}

\begin{center}
{\bf 2A. Non-crossing partitions}
\end{center}

\begin{notation}   \label{def:2.1}
$1^o$ We will use the standard conventions of notation for
non-crossing partitions (as in \cite{S00}, or in Lecture 9 of
\cite{NS06}). For a positive integer $n$, the set of all non-crossing
partitions of $\{ 1, \ldots , n \}$ will be denoted by $NC(n).$
For $\pi \in NC(n)$, the number of blocks of $\pi$ will be denoted
by $| \pi |$. On $NC(n)$ we consider the partial order given
by {\em reversed refinement}: for $\pi , \rho \in NC(n)$, we write
``$\pi \leq \rho$'' to mean that every block of $\rho$ is a union of
blocks of $\pi$. The minimal and maximal element of $( NC(n), \leq )$
are denoted by $0_n$ (the partition of $\{ 1, \ldots , n \}$ into $n$
blocks of 1 element each) and respectively $1_n$ (the partition of
$\{ 1, \ldots , n \}$ into 1 block of $n$ elements).

$2^o$ Every partition $\pi \in NC(n)$ has associated to it a
permutation of $\{ 1, \ldots , n \}$, which is denoted by
$P_{\pi}$, and is defined by the following prescription: for
every block $B = \{ b_1, \ldots , b_m \}$ of $\pi$, with
$b_1 < \cdots < b_m$, one creates a cycle of $P_{\pi}$ by putting
\[
P_{\pi} (b_1) = b_2, \ldots , P_{\pi} (b_{m-1}) = b_m,
P_{\pi} (b_m) = b_1.
\]
Note that in the particular case when $\pi = 0_n$ we have that
$P_{0_n}$ is the identical permutation of $\{ 1, \ldots , n \}$,
while for $\pi = 1_n$ we have that $P_{1_n}$ is the cycle
$1 \mapsto 2 \mapsto \cdots \mapsto n \mapsto 1$.

$3^o$ The {\em Kreweras complementation map} is a special
order-reversing bijection $K: NC(n) \to NC(n)$. In this paper we
will use its description in terms of permutations associated to
non-crossing partitions: for $\pi \in NC(n)$, the Kreweras complement
of $\pi$ is the partition $K( \pi ) \in NC(n)$ uniquely determined
by the fact that its associated permutation is
\begin{equation}  \label{eqn:2.11}
P_{K( \pi )} = P_{\pi}^{-1} P_{1_n}.
\end{equation}
Formula (\ref{eqn:2.11}) can be extended in order to cover the
concept of {\em relative Kreweras complement} of $\pi$ in $\rho$, 
for $\pi, \rho \in NC(n)$ such that $\pi \leq \rho$. This is the 
partition in $NC(n)$, denoted by $K_{\rho} ( \pi )$, uniquely
determined by the fact that the permutation associated to it is
\begin{equation}  \label{eqn:2.12}
P_{K_{\rho} ( \pi )} = P_{\pi}^{-1} P_{\rho}.
\end{equation}
Clearly, the Kreweras complementation map $K$ from (\ref{eqn:2.11})
is the relative complementation with respect to the maximal element
$1_n$ of $NC(n)$.

The formulas (\ref{eqn:2.11}), (\ref{eqn:2.12}) do not follow exactly
the original approach used by Kreweras in \cite{K72}, but are easily
seen to be equivalent to it (see e.g. \cite{NS06}, Exercise 18.25
on p. 301).
\end{notation}

\begin{remark}    \label{rem:2.2}
In this remark we record a few facts about relative Kreweras
complements that will be used later on in the paper.

$1^o$  For $\pi \leq \rho$ in $NC(n)$, the number of blocks of the
relative Kreweras complement $K_{\rho} ( \pi )$ is
\begin{equation}   \label{eqn:2.21}
| K_{\rho} ( \pi )| = n+ | \rho | - | \pi |
\end{equation}
(see e.g. \cite{NS06}, Exercise 18.23 on p. 300).

$2^o$ For a fixed partition $\rho \in NC (n)$, the relative
Kreweras complements $K_{\rho} ( \pi )$ of partitions $\pi \leq \rho$
can be obtained by taking ``separate Kreweras complements'' inside
each block of $\rho$. More precisely, let us write explicitly
$\rho = \{ B_1 , \ldots , B_q \}$. It is easy to see that one has a
natural poset isomorphism
\begin{equation}   \label{eqn:2.22}
\{ \pi \in NC(n) \mid \pi \leq \rho \} \ni \pi \mapsto
( \pi_1, \ldots , \pi_q ) \in
NC( \, |B_1| \, ) \times \cdots \times NC( \, |B_q| \, )
\end{equation}
where for every $1 \leq j \leq q$ the partition
$\pi_j \in NC( |B_j| )$ is obtained by restricting $\pi$ to $B_j$
and by re-denoting the elements of $B_j$, in increasing order, so that
they become $1,2, \ldots , |B_j|$. The above statement about taking
separate complements inside each block of $\rho$ then amounts to the
fact that for $\pi \leq \rho$ in $NC(n)$ we have the implication
\begin{equation}   \label{eqn:2.23}
\Bigl( \, \pi \mapsto
( \pi_1, \ldots , \pi_q ) \, \Bigr) \Rightarrow
\Bigl( \, K_{\rho} ( \pi ) \mapsto
( K( \pi_1 ), \ldots , K( \pi_q) ) \, \Bigr),
\end{equation}
where $K( \pi_1 ), \ldots , K( \pi_q )$ are Kreweras complements
calculated in $NC( |B_1| ), \ldots , NC( |B_q| )$, respectively.
For a discussion of this, see pp. 288-290 in Lecture 18 of
\cite{NS06}.

$3^o$ From part $2^o$ of this remark it is immediate that, for a
fixed $\rho \in NC(n)$, the map $\pi \mapsto K_{\rho} ( \pi )$
is an order-reversing bijection from
$\{ \pi \in NC(n) \mid \pi \leq \rho \}$ onto itself. On the other
hand, it can be shown that
\begin{equation}  \label{eqn:2.24}
\pi \leq \rho_1  \leq \rho_2 \mbox{ in $NC(n)$ } \ \Rightarrow \
K_{\rho_1} ( \pi ) \leq K_{\rho_2} ( \pi )
\end{equation}
(see Exercise 18.25.1 on p. 301 of \cite{NS06}). Hence the
expression ``$K_{\rho} ( \pi )$'' is increasing as a function of
$\rho$, and is decreasing as a function of $\pi$.

The particular case when $\rho_2 = 1_n$ in (\ref{eqn:2.24}) gives
us that
\begin{equation}  \label{eqn:2.25}
\pi \leq \rho \mbox{ in $NC(n)$ } \ \Rightarrow \
K_{\rho} ( \pi ) \leq K ( \pi ).
\end{equation}
Thus for every $\pi \in NC(n)$ it makes sense to define a map
\begin{equation}     \label{eqn:2.26}
\left\{   \begin{array}{rcl}
\{ \rho \in NC(n) \mid  \rho \geq \pi \}
& \rightarrow &
\{ \sigma \in NC(n) \mid \sigma \leq K( \pi ) \}                \\
\rho  & \mapsto &  K_{\rho} ( \pi )  .
\end{array}   \right.
\end{equation}
This map turns out to be bijective, and a helpful fact for
calculating its inverse is that
\begin{equation}   \label{eqn:2.27}
\Bigl( \, \rho \geq \pi, \ K_{\rho} ( \pi ) = \sigma \, \Bigr)
\ \Rightarrow \ K_{K( \pi )} ( \sigma ) = K( \rho )
\end{equation}
(see Lemma 18.9 and Remark 18.10 on p. 291 of \cite{NS06}).
Finally, let us also record here the fact that the whole family of
bijections in (\ref{eqn:2.26}) (corresponding to the various
partitions $\pi \in NC(n)$) can be consolidated into one bijection,
\begin{equation}     \label{eqn:2.28}
\left\{   \begin{array}{rcl}
\{ ( \pi , \rho ) \mid \pi , \rho \in NC(n), \ \pi \leq \rho \}
& \rightarrow &
\{ ( \pi , \sigma ) \mid \pi , \sigma \in NC(n), \
\sigma \leq K( \pi ) \}                                          \\
( \pi , \rho ) & \mapsto & ( \pi, K_{\rho} ( \pi ) ) .
\end{array}   \right.
\end{equation}
\end{remark}

$\ $

\begin{center}
{\bf 2B. Power series}
\end{center}

Here we review a few relevant notations and facts about series in
non-commuting indeterminates, and in particular about multivariable
$R$-transforms.

\begin{notation}    \label{def:2.3}
Let $k$ be a positive integer. As already mentioned in the 
introduction, we use the notation $[k]^*$ for the set of all words 
of finite length over the alphabet $\{ 1, \ldots , k \}$:
\begin{equation}    \label{eqn:2.31}
[k]^* := \cup_{n=0}^{\infty} \, \{ 1, \ldots , k \}^n.
\end{equation}
The length of a word $w \in [k]^*$ will be denoted by $|w|$.
(In (\ref{eqn:2.31}) we followed the standard procedure of also
including into $[k]^*$ a unique word $\phi$ with $| \phi | = 0$.)
\end{notation}

\begin{definition}  \label{def:2.4}
{\em (Series and their coefficients.)}
Let $k$ be a positive integer.

$1^o$ We will use the notation
$\bC_0 \langle \langle z_1, \ldots , z_k \rangle \rangle$
for the set of power series with complex coefficients and with
vanishing constant term in the non-commuting indeterminates
$z_1, \ldots , z_k$. The general form of a series
$f \in \bC_0 \langle \langle z_1, \ldots , z_k \rangle \rangle$ is
thus
\begin{equation}  \label{eqn:2.41}
f( z_1, \ldots , z_k) = \sum_{n=1}^{\infty} \
\sum_{i_1, \ldots , i_n =1}^k  \ \alpha_{(i_1, \ldots , i_n)}
z_{i_1} \cdots z_{i_n} =
\sum_{ \begin{array}{c}
{\scriptstyle w \in [k]^*,}  \\ 
{\scriptstyle |w| \geq 1}
\end{array} } \,  \alpha_w z_w,
\end{equation}
where the coefficients  $\alpha_w$ are from $\bC$ and where, same
as in the introduction, we write in short
$z_w := z_{i_1} \cdots z_{i_n}$ for
$w= (i_1, \ldots , i_n) \in \{ 1, \ldots , k \}^n$, $n \geq 1$.

$2^o$ For every word $w \in [k]^*$ with $|w| \geq 1$ we will
denote by
\[
\cf_{w} :
\bC_0 \langle \langle z_1, \ldots , z_k \rangle \rangle \to \bC
\]
the linear functional which extracts the coefficient of $z_w$
in a series
$f \in \bC_0 \langle \langle z_1, \ldots , z_k \rangle \rangle$.
Thus for $f$ written as in Equation (\ref{eqn:2.41}) we have
$\cf_w (f) = \alpha_w$.

$3^o$ Suppose we are given a positive integer $n$, a word
$w = ( i_1, \ldots , i_n ) \in \{ 1, \ldots , k \}^n$ and a
partition $\pi \in NC(n)$. We define a (generally non-linear)
functional
\[
\cf_{w ; \pi} :
\bC_0 \langle \langle z_1, \ldots , z_k \rangle \rangle \to \bC ,
\]
as follows. For every block $B = \{ b_1, \ldots , b_m \}$ of $\pi$,
with $1 \leq b_1 < \cdots < b_m \leq n$, let us use the notation
\begin{equation}   \label{eqn:2.42}
w \, \vert \, B =
(i_1, \ldots , i_n) \vert B \ := \ (i_{b_1}, \ldots , i_{b_m})
\in \{ 1, \ldots , k \}^m.
\end{equation}
Then we define
\begin{equation}  \label{eqn:2.43}
\cf_{w ; \pi} (f) \ := \
\prod_{B \ \mathrm{block \ of} \ \pi} \ \cf_{w|B}
(f), \ \ \forall \, f \in
\bC_0 \langle \langle z_1, \ldots , z_k \rangle \rangle .
\end{equation}
(For example if $w = (i_1, \ldots , i_5)$ is a word of length 5 and if
$\pi = \{ \{ 1,4,5 \} , \{ 2,3 \} \} \in NC(5)$,
then the above formula comes to
$\cf_{(i_1, i_2, i_3, i_4, i_5) ; \pi } (f)$ =
$\cf_{(i_1, i_4, i_5)} (f) \cdot \cf_{(i_2, i_3)} (f)$,
$f \in \bC_0 \langle \langle z_1, \ldots , z_k \rangle \rangle$.)
\end{definition}

\begin{definition}   \label{def:2.5}
Let $\mu$ be a distribution in $\dalg (k)$ (that is,
$\mu : \bC \langle X_1, \ldots , X_k \rangle \to \bC$ is a linear
functional such that $\mu (1) = 1$). The {\em $R$-transform}
of $\mu$ is the series $R_{\mu}
\in \bC_0 \langle \langle z_1, \ldots , z_k \rangle \rangle$ 
uniquely determined by the requirement that for every $n \geq 1$
and every $1 \leq i_1, \ldots , i_n \leq k$ one has
\begin{equation}   \label{eqn:2.51}
\mu ( X_{i_1} \cdots X_{i_n} ) =
\sum_{\pi \in NC(n)} \cf_{(i_1, \ldots , i_n); \pi} ( R_{\mu} ).
\end{equation}
\end{definition}

\begin{remark}    \label{rem:2.6}
It is easy to see that Equation (\ref{eqn:2.51}) does indeed
determine a unique series in
$\bC_0 \langle \langle z_1, \ldots , z_k \rangle \rangle$. The
coefficients of $R_{\mu}$ are called the free cumulants of $\mu$,
and because of this reason Equation (\ref{eqn:2.51}) is
sometimes referred to as the (free) ``moment--cumulant formula'' 
-- see Lectures 11 and 16 in \cite{NS06}.

Observe that, for $n=1$, Equation (\ref{eqn:2.51}) simply says
that $\cf_{(i)} ( R_{\mu} ) = \mu (X_i)$, $1 \leq i \leq k$. This
explains why for $\mu \in \cG_k$ the $R$-transform $R_{\mu}$ was
displayed in Equation (\ref{eqn:1.6}) of the introduction
by having all its linear coefficients equal to 1.

An important point for the present paper is that the $R$-transform
has a very nice behaviour under the operation $\boxtimes$. This is
recorded in the next proposition.
\end{remark}

\begin{proposition}    \label{prop:2.7}
Let $\mu, \nu$ be distributions in $\dalg (k)$, and let
$w$ be a word in $[k]^*$, with $|w| \geq 1$. Then
\begin{equation}    \label{eqn:2.71}
\cf_w \bigl( \, R_{\mu \boxtimes \nu} \, \bigr)
= \sum_{\pi \in NC(n)} \
\cf_{w ; \pi} \bigl( \, R_{\mu} \, \bigr) \cdot
\cf_{w ; K( \pi )} \bigl( \, R_{\nu} \, \bigr).
\end{equation}
\end{proposition}

For the proof of Proposition \ref{prop:2.7} we refer to
Theorem 14.4 and Proposition 17.2 of \cite{NS06}.

$\ $

\begin{center}
{\bf 2C. Graded connected Hopf algebras}
\end{center}

We will work with graded bialgebras over $\bC$ and we will use the 
standard conventions for notations regarding them (as in the monograph
\cite{S69}, for instance). Our review here does not aim at generality, 
but just covers the specialized Hopf algebras used in 
the present paper.

\begin{notation}   \label{def:2.8}
Let $\cB$ be a graded bialgebra over $\bC$. 

$1^o$ The comultiplication and counit of $\cB$ will be denoted by 
$\Delta$ and respectively $\ee$ (or by $\Delta_{\cB}$ and $\ee_{\cB}$ 
when necessary to distinguish $\cB$ from other graded bialgebras that 
are considered at the same time). The iterations of $\Delta$ will be
denoted as $\Delta^{\ell}$, $\ell \geq 1$. Thus every $\Delta^{\ell}$
is a linear map from $\cB$ to $\cB^{ \otimes \ell}$, where 
$\Delta^1 = \id$ (the identity map from $\cB$ to $\cB$), 
$\Delta^2 = \Delta$, and for $\ell \geq 3$ we put
\begin{equation}  \label{eqn:2.81}
\Delta^{\ell} :=   \bigl( \Delta \otimes
\underbrace{\id \otimes \cdots \otimes \id}_{\ell-2} \, \bigr)
\circ \Delta^{\ell -1}.
\end{equation}  

$2^o$ For every $n \geq 0$, the vector subspace of $\cB$ which 
consists of homogeneous elements of degree $n$ will be denoted by 
$\cB_n$. We thus have a direct sum decomposition 
$\cB = \oplus_{n=0}^{\infty} \cB_n$ where
\[
\left\{  \begin{array}{l}
\cB_0 \ni 1_{\cB} \mbox{ (the unit of $\cB$), }             \\
\cB_m \cdot \cB_n \subseteq \cB_{m+n}, \ \ \forall \, m,n \geq 0,
\end{array} \right.  \ \ \mbox{ and } \ \ 
\left\{  \begin{array}{l}
\ee \vert \, { }_{\cB_n} = 0, \ \ \forall \, n \geq 1,      \\
\Delta ( \cB_n ) \subseteq \oplus_{i=0}^n
\cB_i \otimes \cB_{n-i}, \ \ \forall \, n \geq 0.
\end{array}  \right.
\]
If the space $\cB_0$ of homogeneous elements of degree 0 is 
equal to $\bC 1_{\cB}$ then we say that the graded 
bialgebra $\cB$ is {\em connected}.
\end{notation}

\begin{remark}  \label{rem:2.9}
{\em (The convolution algebra $L( \cB , \cM )$.)}
Let $\cB$ be a graded connected bialgebra, let $\cM$ be a unital 
algebra over $\bC$, and let $L( \cB, \cM )$ denote the vector space
of all linear maps from $\cB$ to $\cM$. For 
$\xi , \eta \in L( \cB , \cM )$ one can define their convolution
product, denoted here simply as ``$\xi \eta$'', by the formula
\begin{equation}  \label{eqn:2.91}
\xi \eta := \mult \circ ( \xi \otimes \eta ) \circ \Delta,
\end{equation}
where $\mult : \cM \otimes \cM \to \cM$ is the linear map given
by multiplication ($\mult (x \otimes y) = xy$ for $x,y \in \cM$).
In other words, (\ref{eqn:2.91}) says that for $b \in \cB$ with
$\Delta (b) = \sum_{i=1}^n b_i ' \otimes b_i ''$ one has
\begin{equation}    \label{eqn:2.92}
( \xi \eta ) (b) := \sum_{i=1}^n \xi (b_i') \eta (b_i'') \in \cM.
\end{equation}
When endowed with its usual vector space structure and with the 
convolution product, $L( \cB , \cM )$ becomes itself a unital
algebra over $\bC$. The unit of $L( \cB , \cM )$ is the
linear map $\cB \ni b \mapsto \ee (b) 1_{\cM}$, which by a slight 
abuse of notation is still denoted as $\ee$ (same notation as for
the counit of $\cB$).

Let us also record here that when one considers a convolution 
product of $\ell \geq 1$ elements 
$\xi_1, \ldots , \xi_{\ell} \in L ( \cB , \cM )$, then
Equation (\ref{eqn:2.91}) becomes
\begin{equation}    \label{eqn:2.93}
\xi_1 \xi_2 \cdots \xi_{\ell} =
\mult_{\ell} \circ ( \xi_1 \otimes \cdots \otimes \xi_{\ell} )
\circ \Delta^{\ell},
\end{equation}
where $\Delta^{\ell} : \cB \to \cB^{\otimes \ell}$ is the $\ell$th 
iteration of the comultiplication and
$\mult_{\ell} : \cM^{ \otimes \ell} \to \cM$ is given by 
multiplication.

For every $\xi \in L( \cB , \cM )$ it makes sense to form 
polynomial expressions in $\xi$, that is, expressions of the form
\[
\sum_{\ell = 0}^n t_{\ell} \xi^{\ell} \in L( \cB , \cM ),
\ \ \mbox{ for $n \geq 0$ and $t_0, t_1, \ldots , t_n \in \bC$, }
\]
where the 
$\xi^{\ell}$ ($0 \leq \ell \leq n$) are convolution powers of $\xi$,
and we make the convention that $\xi^0 := \ee$ (the unit of 
$L( \cB , \cM )$). An important point for the present paper is that
if $\xi$ is such that $\xi ( 1_{\cB} ) = 0$ then it also makes 
sense to define an element 
\begin{equation}   \label{eqn:2.94}
\eta :=
\sum_{\ell = 0}^{\infty} t_{\ell} \xi^{\ell} \in L( \cB , \cM ),
\end{equation}
for an arbitrary infinite sequence $( t_{\ell} )_{\ell \geq 0}$ in 
$\bC$. Indeed, if $\xi ( 1_{\cB} ) = 0$ then by using 
(\ref{eqn:2.93}) and the fact that $\Delta$ respects the grading
one immediately sees that $\xi^{\ell}$ vanishes on $\cB_n$ 
whenever $\ell > n$. Thus $\eta$ from (\ref{eqn:2.94}) can be 
defined as the unique linear map from $\cB$ to $\cM$ which satisfies
\begin{equation}   \label{eqn:2.95}
\eta \, \vert \, { }_{\cB_n} =  \Bigl( \, 
\sum_{\ell = 0}^N t_{\ell} \xi^{\ell} \, \Bigr) \, \vert \,
{ }_{\cB_n}, \ \ \forall \, N \geq n \geq 0.
\end{equation}
\end{remark}

\begin{remark}    \label{rem:2.10}
{\em (The antipode.)}
Here we consider the special case of Remark \ref{rem:2.9} where 
$\cM = \cB$, and we observe that $\cB$ is sure to be a Hopf 
algebra -- this means, by definition, that the identity map 
$\id : \cB \to \cB$ is an invertible element in the convolution 
algebra $L( \cB , \cB )$. The inverse of $\id$ is called the
{\em antipode} of $\cB$ and is denoted by $S$. A reason why
$S$ is sure to exist is that one can introduce it via a 
series expansion as in (\ref{eqn:2.94}) above, which mimics the
geometric series expansion of
$\bigl( \ee - ( \ee - \id ) \bigr)^{-1}$. That is, one can put
\begin{equation}    \label{eqn:2.101}
S := \ee + \sum_{\ell =1}^{\infty} ( \ee - \id )^{\ell}
\in L( \cB , \cB )
\end{equation}
(which makes sense because $( \ee - \id ) ( 1_{\cB} ) = 0$), and
one can then verify that $S$ from (\ref{eqn:2.101}) has indeed 
the property that $S \, \id = \ee = \id \, S$.
See e.g. Propositions 5.2 and 5.3 in \cite{FG-B05}, 
or Theorem 14 in Appendix B of \cite{AEP06}.

A remarkable fact about $S$ (holding for the antipode of any 
Hopf algebra) is that it is a unital anti-homomorphism of $\cB$, 
that is, one has $S( 1_{\cB} ) = 1_{\cB}$ and 
$S(bb') = S(b') S(b), \ \ \forall \, b, b' \in \cB$;
see Proposition 4.0.1 in \cite{S69}.
\end{remark}

\begin{remark}    \label{rem:2.11}
{\em (The group of characters.)}
Here we consider the special case of Remark \ref{rem:2.9} where 
$\cM = \bC$. A unital homomorphism from $\cB$ to $\bC$ is called a 
{\em character}. The set of all characters of $\cB$ will be 
denoted by $\bX ( \cB )$. It is easy to verify that the convolution 
product of two characters is again a character. Moreover, if $\eta$ 
is a character then it is obvious that the functional 
$\eta \circ S \in L( \cB , \bC )$ is a character as well, and it is
easy to verify that 
$\eta \, ( \eta \circ S ) = \ee = ( \eta \circ S ) \, \eta$. 
Hence $\bX ( \cB )$ is a subgroup of the group of 
invertibles of $L( \cB , \bC )$, and is thus referred to as the 
{\em group of characters} of $\cB$. 

The fact observed above, that the inverse of $\eta \in \bX ( \cB )$ 
is $\eta \circ S$, can be combined with Equation (\ref{eqn:2.101}) 
from Remark \ref{rem:2.10} to give us a formula for calculating
inverses of characters via a series expansion,
\begin{equation}    \label{eqn:2.111}
\eta^{-1}  = \ee  + \sum_{\ell =1}^{\infty} ( \ee - \eta )^{\ell} ,
\ \ \mbox{ for $\eta \in \bX ( \cB )$. }
\end{equation}
\end{remark}

\begin{remark}   \label{rem:2.12}
{\em (Exponentials and logarithms for functionals.)}
Let $\cB$ be a graded connected Hopf algebra. If $\xi$ is a 
functional in $L( \cB , \bC )$ such that $\xi ( 1_{\cB} ) = 0$, 
then it makes sense to define its exponential by the familiar
formula
\begin{equation}    \label{eqn:2.121}
\exp \xi =  \sum_{\ell =0}^{\infty} \frac{1}{\ell !} \xi^{\ell},
\end{equation}
where on the right-hand side of (\ref{eqn:2.121}) we use the 
recipe for series reviewed in Remark \ref{rem:2.9}. It is easy
to see that exp maps bijectively the set of functionals 
$\{  \xi \in L( \cB , \bC ) \mid \xi ( 1_{\cB} ) = 0 \}$ onto 
$\{  \eta \in L( \cB , \bC ) \mid \eta ( 1_{\cB} ) = 1 \}$; the 
inverse of this bijection is denoted as ``log'' and can be 
described by using the Taylor series expansion for logarithm:
\begin{equation}    \label{eqn:2.122}
\log \eta = - \sum_{\ell =1}^{\infty} \frac{1}{\ell}
( \ee - \eta )^{\ell}, \ \ \
\mbox{ for $\eta \in L( \cB , \bC )$ with $\eta (1_{\cB}) = 1$.}
\end{equation}
By adjusting the familiar argument for the exponential of a sum of 
two matrices, one finds that
\begin{equation}    \label{eqn:2.213}
\exp ( \xi_1 + \xi_2 )  = \exp ( \xi_1 )  \, \exp ( \xi_2 ), \ \
\forall \, \xi_1, \xi_2 \in L( \cB , \bC ) 
\begin{array}[t]{l}
\mbox{ such that } \xi_1 ( 1_{\cB} ) = \xi_2 ( 1_{\cB} ) = 0  \\
\mbox{ and such that $\xi_1 \xi_2 = \xi_2 \xi_1$.}
\end{array}  
\end{equation}
As a consequence, in the opposite direction of the exp/log bijection
one finds that
\begin{equation}    \label{eqn:2.214}
\log ( \eta_1  \eta_2 )  = \log ( \eta_1 ) + \log ( \eta_2 ), \ \ 
\forall \, \eta_1, \eta_2 \in L( \cB , \bC ) 
\begin{array}[t]{l}
\mbox{ such that } \eta_1 ( 1_{\cB} ) = \eta_2 ( 1_{\cB} ) = 1  \\
\mbox{ and such that $\eta_1 \eta_2 = \eta_2 \eta_1$.}
\end{array}
\end{equation}

Finally, let us record here that one can also adjust to our current 
framework the well-known exp/log correspondence between 
derivations and homomorphisms from Lie algebra theory.
Let $\bX ( \cB )$ be the group of characters of $\cB$ from the 
preceding remark, and on the other hand let us consider the 
set $\bI ( \cB )$ of {\em infinitesimal characters} (sometimes 
also called {\em $\ee$-derivations}) of $\cB$, 
\begin{equation}   \label{eqn:2.125}
\bI ( \cB ) := \{ \xi \in L( \cB , \bC ) \mid
\xi (bb') = \xi (b) \ee (b') + \ee (b) \xi (b'), \ \ 
\forall \, b,b' \in \cB \} .
\end{equation}
Then $\bX ( \cB )$ is contained in the domain of the log map 
discussed above, $\bI ( \cB )$ is contained in the domain of the 
exp map, and the two sets correspond to each other in the sense 
that one has $\log \bigl( \, \bX ( \cB ) \, \bigr) = \bI ( \cB )$.
For a discussion of this fact, see for instance Section 4 of 
the survey paper \cite{FG-B05}.
\end{remark}

$\ $

$\ $

\begin{center}
{\bf\large 3. The Hopf algebra \boldmath{$\cYk$} }
\end{center}
\setcounter{section}{3}
\setcounter{equation}{0}
\setcounter{theorem}{0}

Throughout this section we fix a positive integer $k$. Same as
in the introduction, we use the notation $\cYk$ for the
commutative algebra of polynomials
\[
\cYk := \bC \bigl[ \, Y_w \mid 
        w \in [k]^*, \ |w| \geq 2 \, \bigr] .
\]
In addition to that, we will also use
the following conventions of notation.

\begin{notation}  \label{def:3.1}

$1^o$ For a word $w \in [k]^*$ such that $|w| =1$ (i.e. such that
$w = (i)$ for some $1 \leq i \leq k$) we put
$Y_w := 1$ (the unit of $\cYk$).

$2^o$ Let $w$ be a word in $[k]^*$ with $|w|=n \geq 1$, and let
$\pi = \{ A_1, \ldots , A_q \}$ be a partition in $NC(n)$. We will
denote
\begin{equation}  \label{eqn:3.11}
Y_{w; \pi} := Y_{w_1} \cdots Y_{w_q} \in \cYk,
\end{equation}
where $w_j = w \, \vert \, A_j$ for $1 \leq j \leq q$ (and where
the restriction $w \, \vert \, A$ of the word $w$ to a non-empty
subset $A \subseteq \{ 1, \ldots , n \}$ is defined in the same
way as in Equation (\ref{eqn:2.42}) of Definition \ref{def:2.4}).
\end{notation}

The comultiplication and counit of $\cYk$ are defined as follows.

\begin{definition}  \label{def:3.2}
$1^o$ Let $\Delta : \cYk \to \cYk \otimes \cYk$ be the unital
algebra homomorphism uniquely determined by the requirement
that for every $w \in [k]^*$ with $|w| = n \geq 2$ we have
\begin{equation}  \label{eqn:3.21}
\Delta ( Y_w ) = \sum_{\pi \in NC(n)} \
Y_{w; \pi} \otimes Y_{w; K( \pi )},
\end{equation}
where we use the conventions of notation introduced above (cf.
Equation (\ref{eqn:3.11})), and where $K( \pi )$ is the Kreweras
complement of a partition $\pi \in NC(n)$.

$2^o$ Let $\ee : \cYk \to \bC$ be the unital algebra homomorphism
uniquely determined by the requirement that
\begin{equation}  \label{eqn:3.22}
\ee ( Y_w ) = 0, \ \ \forall \, w \in [k]^*
\mbox{ with $|w| \geq 2$.}
\end{equation}
\end{definition}

We will verify that the maps $\Delta$ and $\ee$ defined above
give indeed a coalgebra (hence bialgebra) structure on $\cYk$.
The only point that is somewhat laborious is the verification of
the coassociativity for $\Delta$, which is covered by the next
two lemmas.

\begin{lemma}   \label{lemma:3.3}
Let $w$ be a word in $[k]^*$ with $|w| =: n \geq 2$. Let $\rho$
be a partition in $NC(n)$, and consider the element
$Y_{w; \rho} \in \cYk$ introduced in Notation \ref{def:3.1}.2.
We have
\begin{equation}  \label{eqn:3.31}
\Delta ( Y_{w; \rho} ) = \sum_{\begin{array}{c}
{\scriptstyle \pi \in NC(n),}  \\
{\scriptstyle \pi \leq \rho}
\end{array}  }  \ Y_{w; \pi} \otimes Y_{w; K_{\rho} ( \pi )},
\end{equation}
where $K_{\rho} ( \pi )$ denotes the relative Kreweras complement
of $\pi$ in $\rho$ (as reviewed in Notation \ref{def:2.1}.3).
\end{lemma}

\begin{proof}
Let us write explicitly $\rho = \{ B_1, \ldots , B_q \}$ and let
us denote $w_j := w \, \vert \, B_j$ for $1 \leq j \leq q$. Then
$Y_{w; \rho} = Y_{w_1} \cdots Y_{w_q}$, hence
\[
\Delta (Y_{w; \rho})
= \Delta (Y_{w_1}) \cdots \Delta (Y_{w_q})
= \prod_{j=1}^q \Bigl( \, \sum_{\pi_j \in NC( |B_j| )} \,
Y_{w_j; \pi_j} \otimes Y_{w_j; K( \pi_j )} \, \Bigr) ,
\]
\begin{equation}    \label{eqn:3.32}
= \sum_{\begin{array}{c}
{\scriptstyle  \pi_1 \in NC( |B_1| ), \ldots }   \\
{\scriptstyle  \ldots , \pi_q \in NC( |B_q| ) }
\end{array} } \
\Bigl( \prod_{j=1}^q Y_{w_j; \pi_j} \Bigr) \otimes
\Bigl( \prod_{j=1}^q Y_{w_j; K( \pi_j )} \, \Bigr) ,
\end{equation}
where the Kreweras complements $K( \pi_j )$ in the latter
expressions are taken in the lattices $NC( |B_j| )$,
$1 \leq j \leq q$.

Now let us consider the bijection (\ref{eqn:2.22}) from
Remark \ref{rem:2.2}.2. It is immediate that if
$\pi \leftrightarrow ( \pi_1, \ldots , \pi_q )$ via this
bijection, then
\[
Y_{w; \pi} = Y_{w_1; \pi_1} \cdots Y_{w_q; \pi_q}.
\]
Moreover, it was noted in the same Remark \ref{rem:2.2}.2 that 
if $\pi \leftrightarrow ( \pi_1, \ldots , \pi_q )$ then
$K_{\rho} ( \pi ) \leftrightarrow ( K( \pi_1 ), \ldots ,
K( \pi_q) )$, hence in this situation we also have that
\[
Y_{w; K_{\rho} ( \pi )} = Y_{ w_1; K( \pi_1 ) } \cdots
Y_{ w_q; K( \pi_q ) }.
\]
Thus when in (\ref{eqn:3.32}) we perform the change of variable
given by the bijection from (\ref{eqn:2.22}), we arrive precisely
to the right-hand side of (\ref{eqn:3.31}), as required.
\end{proof}

\begin{lemma}   \label{lemma:3.4}
The map $\Delta$ from Definition \ref{def:3.2}.1 is coassociative.
\end{lemma}

\begin{proof} It suffices to check that
$( \Delta \otimes \id ) \circ \Delta$ and
$( \id \otimes \Delta ) \circ \Delta$ agree on every generator
$Y_w$, with $w \in [k]^*$ such that $|w| \geq 2$.
So let us fix such a $w$, and let us denote $|w| = n$. On the
one hand we have
\[
\Bigl( \, ( \Delta \otimes \id ) \circ \Delta \, \Bigr) (Y_w)
= ( \Delta \otimes \id )
\Bigl( \, \sum_{\rho \in NC(n)} Y_{w; \rho} \otimes
Y_{w; K( \rho ) } \, \Bigr)
\]
\[
= \sum_{\rho \in NC(n)} \ \Bigl( \,
\sum_{ \begin{array}{c}
{\scriptstyle \pi \in NC(n), } \\
{\scriptstyle \pi \leq \rho }
\end{array}  } \ Y_{w; \pi} \otimes Y_{w; K_{\rho} ( \pi )} \, \Bigr)
\otimes Y_{w; K( \rho ) }  \ \ \mbox{ (by Lemma \ref{lemma:3.3})}
\]
\begin{equation}  \label{eqn:3.41}
= \sum_{ \begin{array}{c}
{\scriptstyle \pi , \rho \in NC(n), } \\
{\scriptstyle \pi \leq \rho }
\end{array}  } \ Y_{w; \pi} \otimes Y_{w; K_{\rho} ( \pi )}
\otimes Y_{ w; K( \rho ) } .
\end{equation}
On the other hand:
\[
\Bigl( \, ( \id \otimes \Delta ) \circ \Delta \, \Bigr) (Y_w)
= ( \id \otimes \Delta )
\Bigl( \, \sum_{\pi \in NC(n)} Y_{w; \pi} \otimes
Y_{w; K( \pi ) } \, \Bigr)
\]
\[
= \sum_{\pi \in NC(n)} Y_{w; \pi} \otimes \Bigl( \,
\sum_{ \begin{array}{c}
{\scriptstyle \sigma \in NC(n), } \\
{\scriptstyle \sigma \leq K( \pi ) }
\end{array}  } \ Y_{w; \sigma} \otimes
Y_{w; K_{K ( \pi )} ( \sigma )} \, \Bigr)
\ \ \mbox{ (by Lemma \ref{lemma:3.3}) }
\]
\begin{equation}  \label{eqn:3.42}
= \sum_{ \begin{array}{c}
{\scriptstyle \pi , \sigma \in NC(n), } \\
{\scriptstyle \sigma \leq K( \pi) }
\end{array}  } \ Y_{w; \pi} \otimes Y_{w; \sigma}
\otimes Y_{w; K_{K ( \pi )} ( \sigma )}.
\end{equation}
Let us now observe that the combinatorial structures
which index the sums from (\ref{eqn:3.41}) and (\ref{eqn:3.42}) are
precisely those appearing in the bijection (\ref{eqn:2.28}) of
Remark \ref{rem:2.2}.3. Moreover, we see (by also taking
Equation (\ref{eqn:2.27}) into account) that the said bijection
(\ref{eqn:2.28}) produces a term-by-term identification of the
sums (\ref{eqn:3.41}) and (\ref{eqn:3.42}); hence
$( \Delta \otimes \id ) \circ \Delta$ and
$( \id \otimes \Delta ) \circ \Delta$ agree on $Y_w$, as
we wanted.
\end{proof}

On $\cYk$ we will also consider a grading, which is defined 
such that every generator $Y_w$ of $\cYk$ gets to be homogeneous
of degree $|w| -1$. More precisely, the homogeneous subspaces 
$\cYk_n$ of $\cYk$ are defined as follows.

\begin{notation}  \label{def:3.5}
For every $n \geq 0$ we denote
\begin{equation}  \label{eqn:3.51}
\cYk_n := \mbox{span} \left\{ Y_{w_1} \cdots Y_{w_q}
\begin{array}{cl}
\vline & q \geq 1, \ w_1, \ldots , w_q \in [k]^* \mbox{ with }  \\
\vline & |w_1|, \ldots , |w_q| \geq 1 \mbox{ and }
|w_1|+ \cdots + |w_q| =n+q
\end{array}  \right\}  .
\end{equation}
\end{notation}

\begin{proposition}  \label{prop:3.6}
With the comultiplication, counit and grading defined above,
$\cYk$ becomes a graded connected Hopf algebra. 
\end{proposition}

\begin{proof}
The coassociativity of the comultiplication $\Delta$ was proved
in Lemma \ref{lemma:3.4}. Let us also verify that $\ee$ and
$\Delta$ satisfy the counit condition. Clearly, it suffices to
check that for every $w \in [k]^*$ with $|w| \geq 2$ we have
\[
\Bigl( \, ( \ee \otimes \id ) \circ \Delta \, \Bigr) (Y_w)
= Y_w =
\Bigl( \, ( \id \otimes \ee ) \circ \Delta \, \Bigr) (Y_w).
\]
For such a word $w$ observe first (directly from how
$\ee$ is defined) that
\[
\ee ( Y_{w; \pi} ) = 0, \ \ \forall \, n \geq 2, \ \forall
\, \pi \in NC(n) \setminus \{ 0_n \}.
\]
But then:
\[
\Bigl( \, ( \ee \otimes \id ) \circ \Delta \, \Bigr) (Y_w)
= ( \ee \otimes \id ) \Bigl( \sum_{\pi \in NC(n)}
Y_{w; \pi} \otimes Y_{w; K( \pi )} \, \Bigr)
\]
\[
= \sum_{\pi \in NC(n)} \ee ( Y_{w; \pi} ) \cdot Y_{w; K( \pi )}
= \ee ( Y_{w; 0_n } ) \cdot Y_{w; 1_n} = Y_{w}.
\]
The verification that
$\Bigl( \, ( \id \otimes \ee ) \circ \Delta \, \Bigr) (Y_w) = Y_w$
is analogous.

We next verify the conditions pertaining to the grading of $\cYk$.
It is immediate that the subspaces $\cYk_n$ ($n \geq 0$) defined
by Equation (\ref{eqn:3.51}) give a direct sum decomposition of
$\cYk$, and satisfy the multiplication relation
\begin{equation}  \label{eqn:3.52}
\cYk_m \cdot \cYk_n \subseteq \cYk_{m+n}, \ \ \forall \,
m,n \geq 0.
\end{equation}
We have to prove that they also satisfy the comultiplication relation
\begin{equation}  \label{eqn:3.53}
\Delta ( \cYk_n ) \subseteq \oplus_{m=0}^n \ \cYk_{m}
\otimes \cYk_{n-m}, \ \ \forall \, n \geq 0.
\end{equation}
By taking into account the specifics of (\ref{eqn:3.51}) and
(\ref{eqn:3.52}), we immediately see it is in fact sufficient
to verify that for every $w \in [k]^*$ with $|w| \geq 2$ we have
\begin{equation}  \label{eqn:3.54}
\Delta ( Y_w ) \in \oplus_{m=0}^{|w|-1} \ \cYk_{m}
\otimes \cYk_{|w|-1-m}.
\end{equation}
And indeed, (\ref{eqn:3.54}) follows from the obvious implication
\[
\left\{  \begin{array}{c}
w \in [k]^* \mbox{ with } |w| = n \geq 2  \\
\pi \in NC(n)
\end{array}  \right\} \Rightarrow
Y_{w; \pi} \otimes Y_{w; K( \pi )} \in
\cYk_{n- | \pi |} \otimes \cYk_{n- | K( \pi ) |},
\]
combined with the fact that for every $\pi \in NC(n)$ we have
$( n- | \pi | ) + ( n- | K( \pi ) | ) = n-1$. (For the latter
fact we invoke the particular case $\rho = 1_n$ of Equation
(\ref{eqn:2.21}) in Remark \ref{rem:2.2}.1.)

Since from (\ref{eqn:3.51}) it is clear that $\cYk_0 = \bC \cdot 1$,
we have thus proved that $\cYk$ is a graded connected bialgebra.
As reviewed in Remark \ref{rem:2.10}, it then automatically follows 
that $\cYk$ is a Hopf algebra.
\end{proof}

We conclude this section by looking at the convolution
of characters of $\cYk$, and by observing the group isomorphism 
stated in Theorem \ref{thm:1.2} from the introduction.

\begin{proposition}  \label{prop:3.7}
For every distribution $\mu \in \cG_k$, let the character
$\chi_{\mu}$ of $\cYk$ be as in Definition \ref{def:1.1}.
Then the map
\begin{equation}   \label{eqn:3.71}
\cG_k \ni \mu \mapsto \chi_{\mu} \in \bX ( \cYk )
\end{equation}
is a group isomorphism, where $\cG_k$ is endowed with the
operation $\boxtimes$ while $\bX ( \cYk )$ is endowed with
the operation of convolution.
\end{proposition}

\begin{proof} It is immediately verified that, for every character
$\xi \in \bX ( \cYk )$, the series
\[
\sum_{i=1}^k z_i + \sum_{ \begin{array}{c}
{\scriptstyle  w \in [k]^*,}   \\
{\scriptstyle  |w| \geq 2}
\end{array}  } \, \xi ( Y_w ) z_w
\]
appears as $R$-transform for a uniquely determined distribution
$\mu \in \cG_k$. This implies that the map (\ref{eqn:3.71}) is a
bijection.

In order to verify that the map (\ref{eqn:3.71}) is a group
homomorphism, we fix two distributions $\mu , \nu \in \cG_k$
for which we prove that
\begin{equation}   \label{eqn:3.72}
\chi_{\mu \boxtimes \nu} = \chi_{\mu} \chi_{\nu} .
\end{equation}
Clearly, it suffices to prove that
the characters appearing on the two sides of (\ref{eqn:3.72})
agree on $Y_w$ for every $w \in [k]^*$ with $|w| \geq 2$. So 
we fix such a $w$, we denote $|w| =:n$, and we compute:
\begin{align*}
\chi_{\mu \boxtimes \nu}  (Y_w)
& = \cf_w ( R_{\mu \boxtimes \nu} ) \ \
\mbox{ (by the definition of $\chi_{\mu \boxtimes \nu}$) }     \\
& = \sum_{\pi \in NC(n)} \cf_{w; \pi} ( R_{\mu} ) \cdot
\cf_{w; K( \pi )} ( R_{\nu} ) \ \
\mbox{ (by Proposition \ref{prop:2.7}) }                       \\
& = \sum_{\pi \in NC(n)} \chi_{\mu} (Y_{w; \pi}) \cdot
\chi_{\nu} (Y_{w; K( \pi )})
\mbox{ (by definition of $\chi_{\mu}, \chi_{\nu}$) }           \\
& = ( \chi_{\mu} \chi_{\nu} ) (Y_w),
\end{align*}
where at the last equality sign of the above calculation we used 
the formula (\ref{eqn:2.92}) for a convolution product, combined
with the formula (\ref{eqn:3.21}) for $\Delta (Y_w)$.
\end{proof}

$\ $

$\ $

\begin{center}
{\bf\large 4. Combinatorial description of
\boldmath{$\log \chi_{\mu}$} }
\end{center}
\setcounter{section}{4}
\setcounter{equation}{0}
\setcounter{theorem}{0}

In this section we continue to fix a positive integer $k$ and
to look at the Hopf algebra $\cYk$. The goal of the section is to
prove Theorem \ref{thm:1.6} from the introduction, which gives
a direct combinatorial description for how the infinitesimal
character $\log \chi_{\mu}$ acts on the generators of $\cYk$.

\begin{definition}   \label{def:4.1}
Let $n$ be a positive integer. A {\em multi-chain} in the lattice
$NC(n)$ is a tuple of the form
\begin{equation}   \label{eqn:4.11}
\cigma = ( \pi_0, \pi_1, \ldots , \pi_{\ell} )
\end{equation}
with $\pi_0, \pi_1, \ldots , \pi_{\ell} \in NC(n)$ such that
$0_n = \pi_0 \leq \pi_1 \leq \cdots \leq \pi_{\ell} = 1_n$. The
positive integer $\ell$ appearing in (\ref{eqn:4.11}) is called
the {\em length} of the multi-chain, and is denoted by
$| \cigma |$. If in (\ref{eqn:4.11}) we have $\pi_{j-1} \neq \pi_j$
for every $1 \leq j \leq \ell$, then we say that $\cigma$ is a
{\em chain} in $NC(n)$.
\end{definition}

In our study of $\cYk$, multi-chains of non-crossing partitions
enter the picture when we consider the iterates
$\Delta^{\ell} : \cYk \to \bigl( \, \cYk \, \bigr)^{\otimes \ell}$
of the comultiplication.

\begin{proposition}  \label{prop:4.2}
Let $w$ be a word in $[k]^*$ with $|w| =: n \geq 2$. Then for every
positive integer $\ell$ we have
\begin{equation} \label{eqn:4.21}
\Delta^{\ell} (Y_w) =
\sum_{\begin{array}{c}
{\scriptstyle \cigma = ( \pi_0, \pi_1, \ldots , \pi_{\ell})}  \\
{\scriptstyle \mathrm{multi-chain} \ \mathrm{in} \ NC(n)}
\end{array} } \
Y_{w; K_{\pi_1} ( \pi_0 )} \otimes
Y_{w; K_{\pi_2} ( \pi_1 )} \otimes \cdots \otimes
Y_{w; K_{\pi_{\ell}} ( \pi_{\ell -1} )}.
\end{equation}
\end{proposition}

\begin{proof}
By induction on $\ell$.
For $\ell =1$ the only multi-chain of length 1 in $NC(n)$ is
$( 0_n, 1_n )$, hence the sum on the right-hand side of
(\ref{eqn:4.21}) has only one term, equal to $Y_w$. Since
$\Delta^1$ is by definition the identity map on $\cYk$, we see
that (\ref{eqn:4.21}) does hold in this case.
For $\ell =2$ we have that $\Delta^2 = \Delta$, while the
right-hand side of (\ref{eqn:4.21}) immediately reduces to the
summation (\ref{eqn:3.21}) used to define $\Delta$ in Section 3.
Hence the formula (\ref{eqn:4.21}) does hold in this case as
well.

Let us now verify the induction step ``$\ell \Rightarrow \ell +1$'',
where $\ell$ is a fixed integer, $\ell \geq 2$. We have that
$\Delta^{\ell + 1} (Y_w) = ( \Delta \otimes
\underbrace{\id \otimes \cdots \otimes \id}_{\ell-1} \, )
( \Delta^{\ell} (Y_w) )$ (cf. Equation (\ref{eqn:2.81}) of 
Notation \ref{def:2.8}). In this expression we replace
$\Delta^{\ell} ( Y_w )$ by using the induction hypothesis; we 
obtain that
\begin{equation}   \label{eqn:4.22}
\Delta^{\ell + 1} (Y_w) =
\sum_{ \begin{array}{c}
{\scriptstyle \cigma = ( \pi_0, \pi_1, \ldots , \pi_{\ell} )}     \\
{\scriptstyle \mathrm{multi-chain} \ \mathrm{in} \ NC(n)}
\end{array} } \
\Delta (Y_{w; \pi_1}) \otimes Y_{w;K_{\pi_2} ( \pi_1 )}
\otimes \cdots \otimes Y_{w;K_{\pi_{\ell}} ( \pi_{\ell -1} )} ,
\end{equation}
where on the right-hand side of (\ref{eqn:4.22}) we also used the
fact that $K_{\pi_1} ( \pi_0 ) = K_{\pi_1} ( 0_n ) = \pi_1$.
In the latter expression we then invoke the explicit formula for
$\Delta ( Y_{w; \pi_1} )$ provided by Lemma \ref{lemma:3.3}. We
obtain that $\Delta^{\ell + 1} (Y_w)$ is equal to
\[
\sum_{ \begin{array}{c}
{\scriptstyle \cigma = ( \pi_0, \pi_1, \ldots , \pi_{\ell} )}     \\
{\scriptstyle \mathrm{multi-chain} \ \mathrm{in} \ NC(n)}
\end{array} } \
\bigl( \, \sum_{ \begin{array}{c}
{\scriptstyle \pi \in NC(n),}     \\
{\scriptstyle \pi \leq \pi_1}
\end{array}  } \ Y_{w; \pi} \otimes
Y_{w; K_{\pi_1} ( \pi )} \, \bigr)
\otimes Y_{w;K_{\pi_2} ( \pi_1 )}
\otimes \cdots \otimes Y_{w;K_{\pi_{\ell}} ( \pi_{\ell -1} )},
\]
hence to
\[
\sum_{ \begin{array}{c}
{\scriptstyle \cigma '=( \pi_0, \pi, \pi_1, \ldots , \pi_{\ell} )} \\
{\scriptstyle \mathrm{multi-chain} \ \mathrm{in} \ NC(n)}
\end{array} } \
Y_{w; K_{\pi} ( \pi_0 ) } \otimes Y_{w; K_{\pi_1} ( \pi )}
\otimes Y_{w;K_{\pi_2} ( \pi_1 )}
\otimes \cdots \otimes Y_{w;K_{\pi_{\ell}} ( \pi_{\ell -1} )},
\]
and this completes the proof of the induction step.
\end{proof}

We now move towards proving the explicit formula for the coefficients
of the $\ls$-transform that was stated in Theorem \ref{thm:1.6}.
We will use ``generalized coefficients'' indexed by multi-chains,
which are defined as follows.

\begin{definition}   \label{def:4.3}
Let $w \in [k]^*$ be a word of length $|w| = n \geq 1$, let
$\cigma = ( \pi_0, \pi_1, \ldots , \pi_{\ell} )$ be a
multi-chain in $NC(n)$, and let $f$ be a series in
$\bC_0 \langle \langle z_1, \ldots , z_k \rangle \rangle$.
The {\em generalized coefficient} of $f$ corresponding
to $w$ and $\cigma$ is
\begin{equation}
\cf_{w}^{ \, ( \cigma )} (f) := \prod_{j=1}^{\ell}
\cf_{ w; K_{\pi_j} ( \pi_{j-1} )} (f).
\end{equation}

\vspace{6pt}

A concrete example: say that $n=4$ and that
$\cigma = ( \pi_0 , \pi_1, \pi_2, \pi_3 )$, where
$\pi_0 = 0_4$, $\pi_3 = 1_4$, and
$\pi_1 = \bigl\{ \, \{ 1,3 \} , \{ 2 \} , \{ 4 \} \, \bigr\}$,
$\pi_2 = \bigl\{ \, \{ 1,3,4 \} , \{ 2 \} \, \bigr\}$.
Then we have
\[
\left\{  \begin{array}{l}
K_{\pi_1} ( \pi_0 ) = \pi_1
= \bigl\{ \, \{ 1,3 \} , \{ 2 \} , \{ 4 \} \, \bigr\},     \\
K_{\pi_2} ( \pi_1 )
= \bigl\{ \, \{ 1 \} , \{ 2 \} , \{ 3,4 \} \, \bigr\} ,    \\
K_{\pi_3} ( \pi_2 ) = K( \pi_2 )
= \bigl\{ \, \{ 1,2 \} , \{ 3 \} , \{ 4 \} \, \bigr\},
\end{array}  \right.
\]
hence for a word $w = (i_1, i_2, i_3, i_4 ) \in [k]^*$ and a series
$f \in \bC_0 \langle \langle z_1, \ldots , z_k \rangle \rangle$
we get that
\[
\cf_{w}^{ \, ( \cigma )} (f) =
\Bigl( \, \cf_{(i_1, i_3)} (f) \, \cf_{(i_2)} (f) \,
\cf_{(i_4)} (f) \, \Bigr) \cdot
\]
\[
\Bigl( \, \cf_{(i_1)} (f) \, \cf_{(i_2)} (f) \,
\cf_{(i_3, i_4)} (f) \, \Bigr) \cdot
\Bigl( \, \cf_{(i_1, i_2)} (f) \, \cf_{(i_3)} (f) \,
\cf_{(i_4)} (f) \, \Bigr) .
\]
Note that if $f$ happens to be an $R$-transform $R_{\mu}$ denoted as
in Equation (\ref{eqn:1.6}) of the introduction, then the above
formula simplifies to
\[
\cf_{w}^{ \, ( \cigma )} (f) = \alpha_{(i_1, i_3)} \cdot
\alpha_{(i_3, i_4)} \cdot \alpha_{(i_1, i_2)}.
\]
\end{definition}

The generalized coefficients defined above come into our
calculations with characters of $\cYk$ via the following lemma.

\begin{lemma}  \label{lemma:4.4}
Let $\mu$ be a distribution in $\cG_k$, let $w$ be a word in 
$[k]^*$ with $|w| = n \geq 2$ and let $\ell$
be a positive integer. If $\ell < n$ then
\begin{equation}   \label{eqn:4.41}
( \chi_{\mu} - \ee )^{\ell} (Y_w) =
\sum_{\begin{array}{c}
{\scriptstyle \cigma \ \mathrm{chain}}  \\
{\scriptstyle \mathrm{in} \ NC(n), \ | \cigma | = \ell}
\end{array} } \ \cf_{w}^{\, ( \cigma )} (R_{\mu}).
\end{equation}
If $\ell \geq n$, then $( \chi_{\mu} - \ee )^{\ell} (Y_w) = 0$.
\end{lemma}

\begin{proof}
We evaluate $( \chi_{\mu} - \ee )^{\ell} (Y_w)$ by using Equation 
(\ref{eqn:2.93}) and then by replacing $\Delta^{\ell} (Y_w)$
from Proposition \ref{prop:4.2}. In this way we arrive to a 
summation over multi-chains,
\begin{equation}   \label{eqn:4.42}
( \chi_{\mu} - \ee )^{\ell} (Y_w) = \sum_{\begin{array}{c}
{\scriptstyle \cigma \ \mathrm{multi-chain} }  \\
{\scriptstyle \mathrm{in} \ NC(n), \ | \cigma | = \ell}
\end{array} } \ \term_{ { }_{\cigma} },
\end{equation}
where the term of the sum indexed by the multi-chain
$\cigma = ( \pi_0, \pi_1, \ldots , \pi_{\ell} )$ is
\begin{equation}   \label{eqn:4.43}
\term_{ { }_{\cigma} } := \prod_{j=1}^{\ell}
( \chi_{\mu} - \ee ) \bigl( Y_{w; K_{\pi_j} ( \pi_{j-1} )} \bigr).
\end{equation}

Note that if in (\ref{eqn:4.43}) there exists an index $j$ such
that $\pi_{j-1} = \pi_{j}$, then it follows that
\[
( \chi_{\mu} - \ee ) \bigl( Y_{w; K_{\pi_j} ( \pi_{j-1} )} \bigr)
= ( \chi_{\mu} - \ee ) ( Y_{w; 0_n} )
= ( \chi_{\mu} - \ee ) ( 1 ) = 1-1 = 0.
\]
This observation shows that $\term_{ { }_{\cigma} }$ vanishes 
whenever $\cigma$ is not a chain.

Since the biggest possible length of a chain in $NC(n)$ is $n-1$,
the above observation implies that for $\ell \geq n$ all the terms
of the sum on the right-hand side of (\ref{eqn:4.42}) are equal
to 0. Hence for $\ell \geq n$ we have that
$( \chi_{\mu} - \ee )^{\ell} (Y_w) = 0$.

Finally, let us assume that $\ell < n$ and let us look at
$\term_{ { }_{\cigma} }$ in the case when
$\cigma = ( \pi_0, \pi_1, \ldots , \pi_{\ell} )$ is
a chain. In this case every relative Kreweras complement
$K_{\pi_j} ( \pi_{j-1} )$ ($1 \leq j \leq \ell$) is different
from $0_n$; this immediately implies that
$\ee \bigl( Y_{w; K_{\pi_j} ( \pi_{j-1} )} \bigr) =0$, hence that
\[
( \chi_{\mu} - \ee ) \bigl( Y_{w; K_{\pi_j} ( \pi_{j-1} )} \bigr)
= \chi_{\mu} \bigl( Y_{w; K_{\pi_j} ( \pi_{j-1} )} \bigr)
= \cf_{ w; K_{\pi_j} ( \pi_{j-1} )} ( R_{\mu} ).
\]
Thus if $\cigma$ is a chain, then the right-hand side of
(\ref{eqn:4.43}) is equal to
$\cf_{w}^{ \, ( \cigma )} ( R_{\mu} )$, and
(\ref{eqn:4.41}) follows.
\end{proof}

As a consequence, we obtain the explicit formula for
$\bigl( \log \chi_{\mu} \bigr) (Y_w)$ that was announced in 
the introduction of the paper.

\begin{proposition}  \label{prop:4.5}
Let $\mu$ be a distribution in $\cG_k$, and let
$w$ be a word in $[k]^*$, with $|w| = n \geq 2$. Then
\begin{equation}    \label{eqn:4.51}
\bigl( \log \chi_{\mu} \bigr) (Y_w) = \sum_{ \begin{array}{c}
{\scriptstyle \cigma \ \mathrm{chain} } \\
{\scriptstyle \mathrm{in} \ NC(n) }
\end{array}  }  \ \frac{ (-1)^{1+| \cigma |} }{| \cigma |}
\cf_{w}^{ \, ( \cigma )} ( R_{\mu} ).
\end{equation}
\end{proposition}

\begin{proof} This follows immediately when we combine the formal 
series expansion which defines $\log \chi_{\mu}$ (Equation 
(\ref{eqn:2.122}) in Remark \ref{rem:2.12}) with the formula
obtained in Lemma \ref{lemma:4.4}.
\end{proof}

\begin{example}    \label{ex:4.6}
Let us see how the expression on the right-hand side of
Equation (\ref{eqn:4.51}) concretely looks, for some small
values of $|w|$.

{\em Case $|w| = 2$.} There exists a unique chain in $NC(2)$,
namely $\cigma = (0_2, 1_2)$. Thus for a word of length 2,
$w = (i_1, i_2)$, the formula (\ref{eqn:4.51}) simply amounts to
\begin{equation}    \label{eqn:4.61}
\cf_{(i_1, i_2)} ( \ls_{\mu} ) = \cf_{(i_1, i_2)} ( R_{\mu} ).
\end{equation}

{\em Case $|w| = 3$.} The lattice $NC(3)$ has 1 chain of length 1,
$\cigma = (0_3, 1_3)$, and has 3 chains of length 2:
$\cigma_i = ( 0_3, \pi_i, 1_3 )$, $1 \leq i \leq 3$, where we put
\[
\pi_1 = \bigl\{ \, \{ 1 \}, \, \{ 2,3 \} \, \bigr\}, \
\pi_2 = \bigl\{ \, \{ 1,3 \}, \, \{ 2 \} \, \bigr\}, \
\pi_3 = \bigl\{ \, \{ 1,2 \}, \, \{ 3 \} \, \bigr\}.
\]
The Kreweras complements of $\pi_1, \pi_2, \pi_3$ are $\pi_2, \pi_3$
and $\pi_1$, respectively. Thus for a word of length 3,
$w = ( i_1, i_2, i_3 )$, we have
\begin{align*}
\cf_{w}^{\, (\cigma_1)} ( R_{\mu} )
& = \cf_{w ; \pi_1} ( R_{\mu} ) \cdot \cf_{w ; \pi_2} ( R_{\mu} ) \\
& = \bigl( \, \cf_{ (i_1, i_2) } ( R_{\mu} )
\cf_{ (i_3) } ( R_{\mu} ) \, \bigr) \cdot
\bigl( \, \cf_{ (i_1, i_3) } ( R_{\mu} )
\cf_{ (i_2) } ( R_{\mu} ) \, \bigr) ,
\end{align*}
and similar formulas hold for $\cf_{w}^{\, (\cigma_2)} ( R_{\mu} )$ 
and $\cf_{w}^{\, (\cigma_3)} ( R_{\mu} )$. By taking into account that
the linear coefficients of $R_{\mu}$ are all equal to 1, we thus get
that
\begin{equation}    \label{eqn:4.62}
\cf_{(i_1, i_2, i_3)} ( \ls_{\mu} )
= \cf_{(i_1, i_2, i_3)} ( R_{\mu} )
- \frac{1}{2} 
\cf_{ (i_1, i_2) } ( R_{\mu} ) \cf_{ (i_1, i_3) } ( R_{\mu} )
\end{equation}
\[
- \frac{1}{2} 
\cf_{ (i_1, i_3) } ( R_{\mu} ) \cf_{ (i_2, i_3) } ( R_{\mu} )
- \frac{1}{2} 
\cf_{ (i_2, i_3) } ( R_{\mu} ) \cf_{ (i_1, i_2) } ( R_{\mu} ).
\]

{\em Case $|w| = 4$.} The lattice $NC(4)$ has 1 chain of length 1,
has 12 chains of length 2 (of the form $( 0_4, \pi , 1_4 )$ with
$\pi$ running in $NC(4) \setminus \{ 0_4, 1_4 \}$), and has 16
chains of length 3. Some direct (but more tedious) calculations
similar to the ones shown above lead to the conclusion that for
a word $w = (i_1, \ldots , i_4) \in \{ 1, \ldots , k \}^4$ one
has
\begin{equation}   \label{eqn:4.63}
\cf_{(i_1, i_2, i_3, i_4)} ( \ls_{\mu} )
= \cf_{(i_1, i_2, i_3, i_4)} ( R_{\mu} )
- \frac{1}{2} \Sigma + \frac{1}{3} \Sigma ' + \frac{1}{6} \Sigma '',
\end{equation}
where each of $\Sigma, \Sigma '$ and $\Sigma ''$ is a sum of
products of coefficients of $R_{\mu}$. (The fraction $1/6$ which
multiplies $\Sigma ''$ appears from cancelations, $1/6$ = $2 \cdot
1/3 - 1/2$.) For the record, we indicate below what these sums are
exactly; in order to save space, the next formula has the
coefficients of $R_{\mu}$ denoted as in Equation (\ref{eqn:1.6}) of
the introduction, by writing ``$\alpha_{ ( \cdots ) }$''
instead of ``$\cf_{ ( \cdots ) } ( R_{\mu} )$''. So, we have:
\begin{equation}  \label{eqn:4.64}
\Sigma = \alpha_{ (i_1, i_2, i_3) } \bigl( \,
\alpha_{ (i_1, i_4) } + \alpha_{ (i_3, i_4) } \, \bigr)
+ \alpha_{ (i_1, i_2, i_4) } \bigl( \,
\alpha_{ (i_3, i_4) } + \alpha_{ (i_2, i_3) } \, \bigr)
\end{equation}
\[
+ \alpha_{ (i_1, i_3, i_4) } \bigl( \,
\alpha_{ (i_2, i_3) } + \alpha_{ (i_1, i_2) } \, \bigr)
+ \alpha_{ (i_2, i_3, i_4) } \bigl( \,
\alpha_{ (i_1, i_2) } + \alpha_{ (i_1, i_4) } \, \bigr),
\]

\begin{equation}  \label{eqn:4.65}
\Sigma ' =
  \alpha_{(i_1, i_2)} \alpha_{(i_1, i_3)} \alpha_{(i_1, i_4)}
+ \alpha_{(i_1, i_2)} \alpha_{(i_2, i_3)} \alpha_{(i_2, i_4)}
+ \alpha_{(i_1, i_3)} \alpha_{(i_2, i_3)} \alpha_{(i_3, i_4)}
\end{equation}
\[
+ \alpha_{(i_1, i_4)} \alpha_{(i_2, i_4)} \alpha_{(i_3, i_4)}
+ \alpha_{(i_1, i_2)} \alpha_{(i_2, i_3)} \alpha_{(i_3, i_4)}
+ \alpha_{(i_1, i_2)} \alpha_{(i_2, i_3)} \alpha_{(i_1, i_4)}
\]
\[
+ \alpha_{(i_1, i_2)} \alpha_{(i_1, i_4)} \alpha_{(i_3, i_4)}
+ \alpha_{(i_1, i_4)} \alpha_{(i_2, i_3)} \alpha_{(i_3, i_4)},
\]
and
\begin{equation}  \label{eqn:4.66}
\Sigma '' =
  \alpha_{(i_1, i_2)} \alpha_{(i_1, i_3)} \alpha_{(i_3, i_4)}
+ \alpha_{(i_1, i_2)} \alpha_{(i_2, i_4)} \alpha_{(i_3, i_4)}
\end{equation}
\[
+ \alpha_{(i_1, i_3)} \alpha_{(i_1, i_4)} \alpha_{(i_2, i_3)}
+ \alpha_{(i_1, i_4)} \alpha_{(i_2, i_3)} \alpha_{(i_2, i_4)}.
\]
\end{example}

$\ $

We conclude this section with the observation that if one removes
the division by $| \cigma |$ on the right-hand side of
(\ref{eqn:4.51}), then what results is a formula for calculating
the $R$-transform of an inverse in the group
$( \cG_k , \boxtimes )$. This happens because (besides its use in
proving Proposition \ref{prop:4.5}) Lemma \ref{lemma:4.4} can also
be used in connection to the inverse formula for characters that 
was reviewed in (\ref{eqn:2.111}) of Remark \ref{rem:2.11}. The
precise statement of this fact goes as follows.

$\ $

\begin{proposition}  \label{prop:4.7}
Let $\mu$ be a distribution in $\cG_k$, and let $\nu$ denote the
inverse of $\mu$ under the operation $\boxtimes$. Then for every
word $w \in [k]^*$ with $|w| = n \geq 2$ we have
\begin{equation}   \label{eqn:4.71}
\cf_w ( R_{\nu} ) = \sum_{ \begin{array}{c}
{\scriptstyle \cigma \ \mathrm{chain} } \\
{\scriptstyle \mathrm{in} \ NC(n) }
\end{array}  }  \ (-1)^{| \cigma |}
\cf_{w}^{ \, ( \cigma )} ( R_{\mu} ).
\end{equation}
\end{proposition}

\begin{proof} We know that
\begin{align*}
\cf_{w} ( R_{\nu} )
& = \chi_{\nu} (Y_w) \ \ \mbox{ (by definition of $\chi_{\nu}$) }   \\
& = \chi_{\mu}^{-1} (Y_w) \ \ \mbox{ (by Theorem \ref{thm:1.2}) }   \\
& = \ee (Y_w) + \sum_{\ell = 1}^{\infty} (-1)^{\ell}
( \chi_{\mu} - \ee )^{\ell} (Y_w) \ \ \mbox{ (by Equation 
                                      (\ref{eqn:2.111})) }         \\
& = \sum_{\ell = 1}^{n-1} (-1)^{\ell} \ \Bigl( \,
\sum_{\begin{array}{c}
{\scriptstyle \cigma \ \mathrm{chain} }  \\
{\scriptstyle \mathrm{in} \ NC(n), \ | \cigma | = \ell}
\end{array} } \ \cf_{w}^{\, ( \cigma )} (R_{\mu}) 
\, \Bigr) ,
\end{align*}
where at the last equality sign we used Lemma \ref{lemma:4.4} and
the fact that $\ee (Y_w) = 0$. The formula (\ref{eqn:4.71})
then immediately follows.
\end{proof}

$\ $

$\ $

\begin{center}
{\bf\large 5. Some basic properties of the LS-transform}
\end{center}
\setcounter{section}{5}
\setcounter{equation}{0}
\setcounter{theorem}{0}

\begin{proposition}   \label{prop:5.1}
Let $( \cA , \varphi )$ be a non-commutative probability space.
Let $a$ be an element of $\cA$ such that $\varphi (a) =1$, and
let $\mu_1 \in \cG_1$ be the distribution of $a$. Moreover, let 
$k$ be a positive integer and let $\mu \in \cG_k$ denote the 
distribution of the $k$-tuple $(a,a, \ldots , a)$. Then 
\begin{equation}  \label{eqn:5.11}
\ls_{\mu} (z_1, \ldots , z_k) = \ls_{\mu_1} (z_1 + \cdots + z_k).
\end{equation}
\end{proposition}

\begin{proof} The $R$-transforms of $\mu$ and of $\mu_1$ are related
by the formula
\begin{equation}  \label{eqn:5.12}
R_{\mu} (z_1, \ldots , z_k) = R_{\mu_1} (z_1+ \cdots + z_k)
\end{equation}
(see Remark 17.13 on p. 280 of \cite{NS06}). This means that 
\begin{equation}  \label{eqn:5.13}
\cf_{(i_1, \ldots , i_n)} (R_{\mu})
= \cf_{ (1, \ldots , 1) } (R_{\mu_1}),
\ \ \forall \, n \geq 1,
\ \forall \, 1 \leq i_1, \ldots , i_n \leq k,
\end{equation}
where on the right-hand side of (\ref{eqn:5.13}) (and same for
(\ref{eqn:5.14})--(\ref{eqn:5.16}) below) there are $n$ repetitions
of 1 in ``$(1, \ldots , 1)$''. From (\ref{eqn:5.13}) one immediately 
sees, by simply following definitions, that
\begin{equation}  \label{eqn:5.14}
\cf_{(i_1, \ldots , i_n); \pi} (R_{\mu})
= \cf_{ (1, \ldots , 1) ; \pi } (R_{\mu_1}),
\ \ \forall \, n \geq 1,
\ \forall \, 1 \leq i_1, \ldots , i_n \leq k,
\ \forall \, \pi \in NC(n),
\end{equation}
and consequently that
\begin{equation}  \label{eqn:5.15}
\cf_{(i_1, \ldots , i_n)}^{ \, ( \cigma )} (R_{\mu})
= \cf_{ (1, \ldots , 1) }^{ \, ( \cigma )} (R_{\mu_1}),
\ \ \forall \, n \geq 1,
\ \forall \, 1 \leq i_1, \ldots , i_n \leq k,
\ \forall \, \cigma,
\end{equation}
where in (\ref{eqn:5.15}) $\cigma$ is an arbitrary
chain of partitions from $NC(n)$.
By using Equation (\ref{eqn:5.15}) and the expression found in
Proposition \ref{prop:4.5} for the coefficients of the
$\ls$-transform, we obtain that
\begin{equation}  \label{eqn:5.16}
\cf_{(i_1, \ldots , i_n)} (\ls_{\mu})
= \cf_{ (1, \ldots , 1) } (\ls_{\mu_1}),
\ \ \forall \, n \geq 2,
\ \forall \, 1 \leq i_1, \ldots , i_n \leq k,
\end{equation}
and (\ref{eqn:5.11}) follows.
\end{proof}

\begin{example}  \label{ex:5.2}
Here is an example of how the $LS$-transform can be used in the
$C^*$-framework. Let $( \cA , \varphi )$ be a tracial \footnote{ The
hypothesis that $\varphi$ is a trace isn't in fact necessary for the
example. The only place where the traciality of $\varphi$ is used is
in order to observe that $(b_1, \ldots ,b_k)$ has the same
distribution as $(a_1 p, \ldots , a_k p)$. This statement is obvious
when $\varphi$ is a trace, but is still true (though less obvious)
without assuming that $\varphi$ is a trace -- one can derive it
solely from the hypothesis that $p$ is free from $\{ a_1, \ldots ,
a_k \}$, via a suitable calculation with free cumulants.}
$C^*$-probability space, let $a_1, \ldots , a_k \in \cA$ be
selfadjoint and let $p \in \cA$ be positive, such that $p$ is freely
independent from $\{ a_1, \ldots , a_k \}$ and such that $\varphi
(p) = \varphi (a_1) = \cdots = \varphi (a_k) = 1$. Consider the
selfadjoint elements $b_1, \ldots , b_k \in \cA$ defined by $b_i =
p^{1/2} a_i p^{1/2}$, $1 \leq i \leq k$. It is immediate that the
$k$-tuple $(b_1, \ldots ,b_k)$ has the same distribution as $(a_1 p,
\ldots , a_k p)$, and hence that
\begin{equation}   \label{eqn:5.21}
\mu_{\vec{b}} = \mu_{\vec{a}} \boxtimes \mu_{\vec{p}},
\end{equation}
where $\mu_{\vec{b}}, \mu_{\vec{a}}, \mu_{\vec{p}}$ denote the
distributions of the $k$-tuples $(b_1, \ldots, b_k)$,
$(a_1, \ldots, a_k)$ and $(p, \ldots , p)$, respectively. Moreover,
it is known that $\mu_{\vec{p}}$ belongs to the centre of
$( \cG_k , \boxtimes )$ (see discussion on pp. 277-279 
of \cite{NS06}). Thus Corollary 1.5 applies to
this situation, and if we also take Proposition \ref{prop:5.1}
into account we arrive to the formula
\begin{equation}   \label{eqn:5.22}
\ls_{\vec{b}}(z_1, \ldots, z_k) =
\ls_{\vec{a}}(z_1, \ldots, z_k) +
\ls_{p}(z_1 + \cdots + z_k),
\end{equation}
where we wrote in short $\ls_{\vec{a}}$ and $\ls_{\vec{b}}$ for the
$\ls$-transforms of $\mu_{\vec{a}}$ and $\mu_{\vec{b}}$, and where
$\ls_p$ stands for the $\ls$-transform of the (1-dimensional)
distribution of the element $p$.
\end{example}

\begin{remark}   \label{rem:5.3}
Let $( \cA , \varphi )$ be a non-commutative probability space,
and let $a_1, \ldots, a_k$ be elements of $\cA$ such that
$\varphi (a_1) = \cdots = \varphi (a_k) =1$. Let $\mu$ be the
distribution of the $k$-tuple $(a_1, \ldots , a_k)$; on the
other hand, for every $1 \leq i \leq k$ let $\mu_i$ denote the
distribution of the element $a_i$ (thus $\mu \in \cG_k$,
whereas $\mu_1, \ldots , \mu_k \in \cG_1$). Then for every
$1 \leq i \leq k$ we have that
\begin{equation}   \label{eqn:5.31}
\ls_{\mu} ( 0, \ldots , 0, z_i , 0, \ldots , 0 )
= \ls_{\mu_i} (z_i).
\end{equation}
The easy proof of (\ref{eqn:5.31}) (which is very similar to the
argument shown in the proof of Proposition \ref{prop:5.1}) is left
to the reader.

We next look at the important special case when (in the same
notations as above) $a_1, \ldots , a_k$ form a freely independent
family. Let us recall here a fundamental fact from the theory of
the multivariable $R$-transform: the free independence of
$a_1, \ldots , a_k$ is equivalent to the vanishing of the mixed
coefficients of $R_{\mu}$, hence to the fact that the series
$R_{\mu} (z_1, \ldots , z_k)$ separates the variables (see
Lecture 16 in \cite{NS06}). The next proposition shows that this
feature of the $R$-transform is passed to the $\ls$-transform.
\end{remark}

\begin{proposition}   \label{prop:5.4}
Let $( \cA , \varphi )$, the elements $a_1, \ldots, a_k \in \cA$,
and the distributions  $\mu \in \cG_k$ and
$\mu_1, \ldots , \mu_k \in \cG_1$ be as in the preceding remark.
The following statements are equivalent:

(1) The elements $a_1, \ldots , a_k$ form a freely independent
family in $( \cA , \varphi )$.

(2) The $\ls$-transform of $\mu$ separates the variables. That is,
there exist series $u_1, \ldots , u_k \in \bC [[z]]$ such that
\begin{equation}  \label{eqn:5.41}
\ls_{\mu} (z_1, \ldots , z_k) = u_1 (z_1) + \cdots
+ u_k (z_k).
\end{equation}

(3) One has
\begin{equation}  \label{eqn:5.42}
\ls_{\mu} (z_1, \ldots , z_k) = \ls_{\mu_1} (z_1) + \cdots
+ \ls_{\mu_k} (z_k).
\end{equation}
\end{proposition}

\begin{proof} The implication {\em (3) $\Rightarrow$ (2)} is
trivial, while {\em (2) $\Rightarrow$ (3)} follows immediately from
Remark \ref{rem:5.3}. The bulk of the proof will be devoted to
the equivalence of {\em (1)} and {\em (2)}.

\vspace{6pt}

{\em ``(1) $\Rightarrow$ (2)''.}
We have to show that all the mixed
coefficients of $\ls_{\mu}$ vanish. So let us fix $n \geq 2$
and $1 \leq i_1, \ldots , i_n \leq k$ such that
$\cf_{(i_1, \ldots , i_n)} ( \ls_{\mu} ) \neq 0$; our goal is
to prove that $i_1 = i_2 = \cdots = i_n$.

$ $From the formula for
$\cf_{(i_1, \ldots , i_n)} ( \ls_{\mu} )$ found in Proposition
\ref{prop:4.5} we infer that there must exist
chains $\cigma$ of partitions in $NC(n)$ such that
$\cf_{(i_1, \ldots , i_n)}^{ \, ( \cigma )} ( R_{\mu} ) \neq 0$.
Let us fix such a chain,
$\cigma = ( \pi_0, \pi_1, \ldots , \pi_{\ell} )$. We have that
\[
0 \neq \cf_{(i_1, \ldots , i_n)}^{ \, ( \cigma )} ( R_{\mu} )
= \prod_{j=1}^k
\cf_{(i_1, \ldots , i_n); K_{\pi_j} ( \pi_{j-1} ) } (R_{\mu}),
\]
hence that
\begin{equation}  \label{eqn:5.43}
\cf_{(i_1, \ldots , i_n); K_{\pi_j} ( \pi_{j-1} ) } (R_{\mu}) \neq 0,
\ \ \forall \, 1 \leq j \leq \ell.
\end{equation}

At this point it is convenient to encode our $n$-tuple
$(i_1, \ldots , i_n)$ as a function
$I: \{ 1, \ldots , n \} \to \{ 1, \ldots , k \}$, where
$I(m) := i_m$ for $1 \leq m \leq n$. Then (\ref{eqn:5.43}) can be
combined with the fact (equivalent to the hypothesis in {\em (1)})
that the mixed coefficients of $R_{\mu}$ are all equal to 0, in
order to infer that the function
$I$ is constant along every block of the partition
$K_{\pi_j} ( \pi_{j-1} )$, for every $1 \leq j \leq \ell$.
In view of how the permutation associated to a non-crossing
partition is defined (see Notation \ref{def:2.1}.2), the latter
fact is equivalent to
\begin{equation}   \label{eqn:5.45}
I \circ P_{ K_{\pi_j} ( \pi_{j-1} ) } = I,
\ \ \forall \, 1 \leq j \leq \ell .
\end{equation}

Now, the formula for the permutation associated to a relative
Kreweras complement (Equation (\ref{eqn:2.12}) in Notation
\ref{def:2.1}.3) gives us that
\begin{equation}   \label{eqn:5.44}
P_{ K_{\pi_j} ( \pi_{j-1} ) } =
P_{ \pi_{j-1} }^{-1} \  P_{ \pi_j }, \ \ \forall \,
1 \leq j \leq \ell .
\end{equation}
When we put together (\ref{eqn:5.45}) and (\ref{eqn:5.44}),
we find that
\begin{align*}
I
& = I \circ \Bigl( \,
P_{ K_{\pi_1} ( \pi_0 ) }
P_{ K_{\pi_2} ( \pi_1 ) } \cdots
P_{ K_{ \pi_{\ell} } ( \pi_{\ell -1} ) } \, \Bigr)    \\
& = I \circ \Bigl( \,
\bigl( P_{ \pi_0 }^{-1} P_{\pi_1} \bigr)
\bigl( P_{ \pi_1 }^{-1} P_{\pi_2} \bigr)
\bigl( P_{ \pi_{\ell -1} }^{-1} P_{\pi_{\ell} } \bigr) \, \Bigr)  \\
& = I \circ \bigl( P_{ \pi_0 }^{-1} P_{\pi_{\ell} } \bigr) .
\end{align*}
But $\pi_0 = 0_n$ and $\pi_{\ell} = 1_n$, so the permutation
$P_{ \pi_0 }^{-1} P_{\pi_{\ell} }$ is just the cycle
$1 \mapsto 2 \mapsto \cdots \mapsto n \mapsto 1$. The fact
that $I = I \circ \bigl( P_{ \pi_0 }^{-1} P_{\pi_{\ell} } \bigr)$
then immediately implies that $I$ is a constant function, hence
that $i_1 = i_2 = \cdots = i_n$, as we wanted.

\vspace{6pt}

{\em ``(2) $\Rightarrow$ (1)''.}
We will verify the free independence of $a_1, \ldots , a_k$ by
checking that all the mixed coefficients of $R_{\mu}$ vanish:
\begin{equation}   \label{eqn:5.46}
\left\{  \begin{array}{c}
\cf_{(i_1, \ldots , i_n)} (R_{\mu}) = 0 \mbox{ for every
$n \geq 2$ and every $1 \leq i_1, \ldots , i_n \leq k$}     \\
\mbox{ for which it is not true that $i_1 = i_2 = \cdots = i_n$.}
\end{array}  \right.
\end{equation}
The verification of (\ref{eqn:5.46}) will be made by induction on $n$.

For the base case $n=2$ we just have to note that
$\cf_{(i_1,i_2)} (R_{\mu}) =
\cf_{(i_1,i_2)} ( \ls_{\mu} ) =0$ for every $1 \leq i_1, i_2 \leq k$
such that $i_1 \neq i_2$ (where we used Equation (\ref{eqn:4.61})
from Example \ref{ex:4.6}, and the hypothesis that the mixed
coefficients of $\ls_{\mu}$ are all equal to 0).

The induction step: we assume that (\ref{eqn:5.46}) was proved for
$2,3, \ldots , n-1$, and we prove that it also holds for $n$, where
$n \geq 3$. Fix some indices $1 \leq i_1, \ldots , i_n \leq k$ for
which it is not true that $i_1 = i_2 = \cdots = i_n$. Observe that
Equation (\ref{eqn:4.51}) from Proposition \ref{prop:4.5} can be
re-written in the form
\begin{equation}  \label{eqn:5.47}
\cf_{(i_1, \ldots , i_n)} (R_{\mu})
= \cf_{(i_1, \ldots , i_n)} ( \ls_{\mu})
+ \sum_{\begin{array}{c}
{\scriptstyle \cigma \ \mathrm{chain} \ \mathrm{in} NC(n)}   \\
{\scriptstyle \mathrm{with} \ | \cigma | \geq 2}
\end{array} } \
\frac{ (-1)^{| \cigma |} }{ | \cigma | }
\cf_{(i_1, \ldots , i_n)}^{ \, ( \cigma )} (R_{\mu}).
\end{equation}
Now, if $\cigma$ is a chain of partitions in $NC(n)$ with
$| \cigma | \geq 2$ then we get that
$\cf_{(i_1, \ldots , i_n)}^{ \, ( \cigma )} (R_{\mu}) = 0$ by an
argument very similar to the one used in the proof of
{\em (1) $\Rightarrow$ (2)} (where this time we also use the
induction hypothesis, and the fact the explicit writing of
$\cf_{(i_1, \ldots , i_n)}^{ \, ( \cigma )} (R_{\mu})$ only
involves coefficients of order $\leq (n-1)$ of $R_{\mu}$).
Since we also know that
$\cf_{(i_1, \ldots , i_n)} ( \ls_{\mu}) =0$ (by the hypothesis
that the mixed coefficients of $\ls_{\mu}$ vanish), Equation
(\ref{eqn:5.47}) gives us that
$\cf_{(i_1, \ldots , i_n)} (R_{\mu}) = 0$, and this completes our
inductive argument.
\end{proof}

$\ $

$\ $

\begin{center}
{\bf\large 6. Case of one variable: LS and log S}
\end{center}
\setcounter{section}{6}
\setcounter{equation}{0}
\setcounter{theorem}{0}

This section is devoted to the one-variable framework: we
discuss connections that $( \cG_1 , \boxtimes )$ and $\cYone$
have with symmetric functions, and we prove Theorem \ref{thm:1.8}.

\begin{remark}   \label{rem:6.1}
In this remark we make an update of our notations for series and we
give a more precise review of what is the $S$-transform.

In the one-variable framework our notations for series are
simplified by the fact that we only have to deal with the
commutative algebra $\bC \, [[z]]$. Observe for instance that if
$\mu$ is a distribution in $\cG_1$, then the notation
used in Equation (\ref{eqn:1.6}) of the introduction
for the $R$-transform $R_{\mu}$ now simplifies to just
\begin{equation}  \label{eqn:6.12}
R_{\mu} (z) = z + \sum_{n=2}^{\infty} \alpha_n z^n.
\end{equation}

In what follows we will work with two kinds of ``inverses'' for
series in $\bC \, [[z]]$. On the one hand we have the inverse under
multiplication, which is defined for a series $f(z) =
\sum_{n=0}^{\infty} t_n z^n$ with $t_0 \neq 0$, and will
be denoted as ``$1/f$''. On the other hand we have the inverse under
composition, which is defined for a series
$f(z) = \sum_{n=1}^{\infty} t_n z^n$ with $t_1 \neq 0$, and will
be denoted as ``$f^{\langle -1 \rangle}$''.

The definition of the $S$-transform for a distribution $\mu \in
\cG_1$ goes as follows: one first considers the moment series of
$\mu$,
\begin{equation}   \label{eqn:6.13}
M_{\mu} (z) := z + \sum_{n=2}^{\infty} \mu ( X^n ) z^n 
\in \, \bC \, [[z]],
\end{equation}
then puts
\begin{equation}    \label{eqn:6.14}
S_{\mu} (z) := \frac{1+z}{z} M_{\mu}^{\langle -1 \rangle} (z).
\end{equation}

In this paper we prefer to write $S_{\mu}$ in terms of the
$R$-transform $R_{\mu}$; this is done by a similarly looking
formula,
\begin{equation}  \label{eqn:6.15}
S_{\mu} (z) = \frac{1}{z} R_{\mu}^{\langle -1 \rangle} (z)
\end{equation}
(see \cite{NS06}, Remark 16.18 on p. 270).

It is clear that for $\mu \in \cG_1$, the $S$-transform $S_{\mu} \in
\bC \, [[z]]$ has constant term equal to 1. Hence one can consider
the inverse under multiplication $1/S_{\mu}$, and this is another
series with constant term equal to 1:
\begin{equation}  \label{eqn:6.16}
1/S_{\mu} (z) = 1 + \sum_{n=1}^{\infty} \gamma_n z^n,
\end{equation}
for some coefficients $(\gamma_n)_{n=1}^{\infty}$ in $\bC$. The
series $1/S_{\mu}$ goes sometimes under the name of ``$T$-transform
of $\mu$'' (see \cite{D07}). We will use an explicit formula for how
the coefficients of $R_{\mu}$ are expressed in terms of those of
$1/S_{\mu}$, which is presented in the next lemma. The proof of the
lemma relies on a functional equation that is used in the study of
the $R$-transform and is described as follows. Let $\sns$ and $\tns$
be two sequences of complex numbers such that the $s_n$'s are
expressed in terms of the $t_n$'s by the formula
\begin{equation}  \label{eqn:6.17}
s_n = \sum_{ \begin{array}{c}
{\scriptstyle \pi = \{ A_1, \ldots , A_q \} } \\
{\scriptstyle \mathrm{in} \ NC(n) }
\end{array} } \ t_{|A_1|} \cdots t_{|A_q|}, \ \
\forall \, n \geq 1.
\end{equation}
Then the series $f(z) := \sum_{n=1}^{\infty} t_n z^n$ and
$g(z) := \sum_{n=1}^{\infty} s_n z^n$ satisfy the functional equation
\begin{equation}  \label{eqn:6.18}
g(z) = f \bigl( \, z(1+g(z)) \, \bigr).
\end{equation}
Conversely, if the functional equation (\ref{eqn:6.18}) is
satisfied then it follows that the relations (\ref{eqn:6.17}) hold.
For the proof of these facts, see Theorem 16.15 in \cite{NS06}.
\end{remark}

\begin{lemma}  \label{lemma:6.2}
Let $\mu$ be in $\cG_1$, and consider the series $R_{\mu}$,
$1/S_{\mu} \in \bC \, [[z]]$, with coefficients denoted as in
Equations (\ref{eqn:6.12}) and (\ref{eqn:6.16}), respectively. Then
for every $n \geq 2$ we have that
\begin{equation}  \label{eqn:6.21}
\alpha_n = \sum_{\begin{array}{c}
{\scriptstyle \pi = \{ A_1, \ldots , A_q \} } \\
{\scriptstyle \mathrm{in} \ NC(n-1) }
\end{array} } \ \gamma_{|A_1|} \cdots \gamma_{|A_q|}.
\end{equation}
\end{lemma}

\begin{proof} Consider the series
\begin{equation}  \label{eqn:6.22}
g(z) := \sum_{n=1}^{\infty} \alpha_{n+1} z^n
= \frac{1}{z} R_{\mu} (z) - 1
\end{equation}
and
\begin{equation}  \label{eqn:6.23}
f(z) := \sum_{n=1}^{\infty} \gamma_n z^n = 1/S_{\mu} (z) - 1.
\end{equation}
The required formula (\ref{eqn:6.21}) says that the coefficients of
$g$ are obtained from those of $f$ by summations over non-crossing
partitions, exactly in the way described in Equation
(\ref{eqn:6.17}). Thus we may, equivalently, prove that $f$ and $g$
satisfy the functional equation (\ref{eqn:6.18}).

We will obtain (\ref{eqn:6.18}) by starting from the formula
(\ref{eqn:6.15}) which gives $S_{\mu}$ in terms of $R_{\mu}$. When
explaining how this done, it is convenient to give a name to the
identity series, say that we call it ``$\id$'':
\[
\id \in \bC \, [[z]], \ \ \id (z)=z.
\]
Equation (\ref{eqn:6.15}) says that $\id \cdot S_{\mu} =
R_{\mu}^{\langle -1 \rangle}$, and this can be processed as follows:
\begin{align*}
\id \cdot S_{\mu} = R_{\mu}^{\langle -1 \rangle}
& \Rightarrow ( \id \cdot S_{\mu} ) \circ  R_{\mu} = \id        \\
& \Rightarrow ( \id \circ  R_{\mu} ) \cdot
              ( S_{\mu} \circ  R_{\mu} ) = \id                  \\
& \Rightarrow R_{\mu}   \cdot ( S_{\mu} \circ  R_{\mu} ) = \id  \\
& \Rightarrow R_{\mu} = \id \cdot
              \Bigl( \, 1/( S_{\mu} \circ  R_{\mu} ) \, \Bigr)  \\
& \Rightarrow R_{\mu} = \id \cdot
   \Bigl( \, \bigl( 1/S_{\mu} \bigr) \circ  R_{\mu} ) \, \Bigr).
\end{align*}
By using the latter equality and the relation between $R_{\mu} (z)$
and $g(z)$ we obtain that $g+1$ = $(1/S_{\mu}) \circ R_{\mu}$, or
equivalently, that
\begin{equation}  \label{eqn:6.24}
g = \Bigl( \, (1/S_{\mu}) -1 \, \Bigr) \circ R_{\mu}.
\end{equation}
But (\ref{eqn:6.24}) implies (\ref{eqn:6.18}), since
$(1/S_{\mu}) -1 = f$ and $R_{\mu} (z) = z(1+g(z))$.
\end{proof}

\begin{remark}  \label{rem:6.3}
In this remark we collect a few basic facts about the algebra $\sym$
of symmetric functions. We will deal with three of the commonly used
sequences of generators for $\sym$:

\[
\left\{   \begin{array}{l}
\mbox{ the sequence $\ens$ of elementary symmetric functions, }   \\
                                                                  \\
\mbox{ the sequence $\hns$ of complete homogeneous symmetric 
       functions, and }                                           \\
                                                                  \\
\mbox{ the sequence $\pns$ of power sum symmetric functions. }
\end{array}   \right.
\]
It is well-known that $\sym$ can be identified to the
commutative polynomial algebra generated by either of these three
sequences (see Theorem 7.4.4 and Corollaries 7.6.2, 7.7.2 of
\cite{S99}, noting the slight difference that here we use $\bC$ as
field of scalars). Moreover, one has very nice formulas for passing
from one sequence of generators to another, which are best expressed
in terms of some suitably defined series in $\sym \, [[z]]$. For the
present paper it is convenient to record these formulas as follows:
consider the generating function for the $\ens$,
\begin{equation}    \label{eqn:6.31}
\uE (z) := 1 + \sum_{n=1}^{\infty} e_n z^n,
\end{equation}
and the slightly modified generating functions for $\hns$ and $\pns$
defined by
\begin{equation}    \label{eqn:6.32}
\uH (z) := 1 + \sum_{n=1}^{\infty} (-1)^n h_n z^n,
\end{equation}
and
\begin{equation}    \label{eqn:6.33}
\uP (z) := \sum_{n=1}^{\infty} \frac{(-1)^{n+1}}{n} p_n z^n.
\end{equation}
Then the formulas connecting the sequences $\ens$ and $\hns$ can be
consolidated into the short statement that
\begin{equation}    \label{eqn:6.34}
\uH  = 1/ \uE .
\end{equation}
Equation (\ref{eqn:6.34}) is an equality holding in the algebra of
series $\sym \, [[z]]$, where we use conventions of notation
analogous to those presented for $\bC \, [[z]]$ in Remark
\ref{rem:6.1}. Similarly, the formulas connecting $\ens$ and $\pns$
are consolidated into the equation
\begin{equation}    \label{eqn:6.35}
\uP =  \log \uE ,
\end{equation}
where $\log \uE \in \sym \, [[z]]$ is defined in the standard way
(for series), by using the Taylor series of the logarithm:
\begin{equation}    \label{eqn:6.36}
\log \uE := - \sum_{n=1}^{\infty} \frac{1}{n} (1- \uE )^n.
\end{equation}
For the proofs of these facts, we refer to Sections 7.6 and 7.7 of
\cite{S99}, or to pp. 12-16 in Chapter I.2 of \cite{M79}.

Since in this paper we are dealing with the Hopf algebra structure
of $\sym$, let us also record the formulas for the comultiplication
of the three basic sequences of generators discussed above. For
$\ens$ and $\hns$ one has that
\begin{equation}  \label{eqn:6.37}
\left\{  \begin{array}{lll}
\Delta ( e_n ) & = & \sum_{i=0}^n e_i \otimes e_{n-i},  \\
               &   &                                    \\
\Delta ( h_n ) & = & \sum_{i=0}^n h_i \otimes h_{n-i},
\end{array}  \ \ \forall \, n \geq 1  \right.
\end{equation}
(where we make the convention that $e_0 = h_0 := 1$). On the other
hand, every power sum symmetric function $p_n$ is a primitive
element of $\sym$, in the sense that one has
\begin{equation}  \label{eqn:6.38}
\Delta ( p_n ) = p_n \otimes 1 \, + \, 1 \otimes p_n, \ \ \forall \,
n \geq 1 .
\end{equation}
It is worth keeping in mind that conversely, every primitive 
homogeneous element of $\sym$ must be a scalar multiple of some
$p_n$ -- see \cite{Z81}, Proposition 3.15 on p. 42 and Proposition 
5.3 on p. 75.
\end{remark}

\vspace{10pt}

We now return to the space of distributions $\cG_1$. As already
mentioned in the introduction, we will use the $S$-transform in
order to find a group isomorphism $\mu \mapsto \theta_{\mu}$ from $(
\cG_1 , \boxtimes )$ onto the group $\bX ( \sym )$ of characters of
$\sym$. The definition of $\theta_{\mu}$ is given next. For this
definition recall that, since $\sym$ can be viewed as the
commutative polynomial algebra generated by $\hns$, a character of
$\sym$ can be defined (and is completely determined) by prescribing
at will how we want it to act on the $h_n$'s.

\begin{definition}  \label{def:6.4}
Let $\mu$ be a distribution in $\cG_1$, and consider its
$S$-transform,
\[
S_{\mu} (z) =: 1 + \sum_{n=1}^{\infty} \beta_n z^n \in \bC \, [[z]].
\]
We denote by $\theta_{\mu}$ the character in $\bX ( \sym )$ which is
uniquely determined by the requirement that
\begin{equation}   \label{eqn:6.41}
\theta_{\mu} (h_n) = (-1)^n \beta_n, \ \ \forall \, n \geq 1.
\end{equation}
\end{definition}

\begin{proposition}   \label{prop:6.5}
The map $\mu \mapsto \theta_{\mu}$ is a group isomorphism from $(
\cG_1 , \boxtimes )$ onto the group $\bX ( \sym )$ of characters of
$\sym$, endowed with the operation of convolution.
\end{proposition}

\begin{proof}
It is immediate that the map $\mu \mapsto \theta_{\mu}$ is a
bijection from $\cG_1$ onto $\bX ( \sym )$: it is injective because
a distribution $\mu \in \cG_1$ is completely determined by its
$S$-transform, and it is surjective because every series in $\bC \,
[[z]]$ which has constant term equal to 1 appears as $S_{\mu}$ for
some $\mu \in \cG_1$. So the main point of the proof is to verify
the group homomorphism property, i.e. that
\begin{equation}   \label{eqn:6.51}
\theta_{\mu \boxtimes \nu}  = \theta_{\mu}  \theta_{\nu},
\ \ \forall \, \mu , \nu \in \cG_1.
\end{equation}
We fix $\mu, \nu$ for which we will prove that (\ref{eqn:6.51})
holds. Since both $\theta_{\mu \boxtimes \nu}$ and
$\theta_{\mu}  \theta_{\nu}$ are in $\bX ( \sym )$, it will suffice 
to prove that they agree on the set of generators $\hns$ of $\sym$.

In view of how convolution of characters is defined, and by taking
into account the formula (\ref{eqn:6.37}) for the comultiplication
of $h_n$, we see that:
\begin{equation}  \label{eqn:6.52}
\bigl( \, \theta_{\mu} \theta_{\nu} \, \bigr) (h_n) = \sum_{i=0}^n
\theta_{\mu} (h_i) \cdot \theta_{\nu} (h_{n-i}), \ \ \forall \, n
\geq 1.
\end{equation}
Let us then consider the $S$-transforms of $\mu$ and of $\nu$, and
let us write explicitly their coefficients:
\[
S_{\mu} (z) =: 1 + \sum_{n=1}^{\infty} \beta_n '  z^n, \ \ S_{\nu}
(z) =: 1 + \sum_{n=1}^{\infty} \beta_n '' z^n.
\]
By invoking the definition of $\theta_{\mu}$ and $\theta_{\nu}$ on
the right-hand side of (\ref{eqn:6.52}) we find that for every $n
\geq 1$ we have:
\begin{align*}
\bigl( \, \theta_{\mu} \theta_{\nu} \, \bigr) (h_n) & = \sum_{i=0}^n
\bigl( \, (-1)^i \beta_i ' \, \bigr) \cdot
\bigl( \, (-1)^{n-i} \beta_{n-i} '' \, \bigr)                      \\
& = (-1)^n \cdot \mbox{ (coeff. of $z^n$ in $S_{\mu} S_{\nu}$)}    \\
& = (-1)^n \cdot \mbox{ (coeff. of $z^n$ in $S_{\mu \boxtimes \nu}$);}
\end{align*}
at the last equality sign we used the multiplicativity
property of the $S$-transform, which says precisely that
$S_{\mu \boxtimes \nu} = S_{\mu} \cdot S_{\nu}$. Clearly, the above
sequence of equalities has ended with
$\theta_{\mu \boxtimes \nu} (h_n)$, and this concludes the proof.
\end{proof}

\begin{remark}   \label{rem:6.6}
The reader may have noticed that the proof of Proposition 
\ref{prop:6.5} would work fine (and would look simpler)
without the sign ``$(-1)^n$'' in front of $h_n$. 
Also, with or without signs, the proof would work equally 
well with $e_n$'s instead of the $h_n$'s, since the $e_n$'s
comultiply in exactly the same way as the $h_n$'s in Equation
(\ref{eqn:6.37}).

Hence there were in fact several possible alternatives for the
definition given above to $\theta_{\mu}$. The justification for
the choice selected in Definition \ref{def:6.4} is that it gives
an underlying homomorphism $\cYone \to \sym$ which matches the
``combinatorial Hopf algebra'' structure from \cite{ABS06}. This 
is discussed in more detail in Remark \ref{rem:7.7} and
Proposition \ref{prop:7.8} below. Right now let us only insist on 
the point that it does indeed make sense to look for an ``underlying 
homomorphism'' in connection to the definition of $\theta_{\mu}$;
this is because we have found 
$( \cG_1, \boxtimes )$ to be isomorphic both to $\bX ( \cYone )$ 
(via $\mu \mapsto \chi_{\mu}$) and to $\bX ( \sym )$ 
(via $\mu \mapsto \theta_{\mu}$), and it is natural
to try to connect these two isomorphisms via a homomorphism 
$\Phi : \cYone \to \sym$. The homomorphism $\Phi$ will be 
introduced in Definition \ref{def:6.8} below. Before introducing
it, we would like to make some comments about $\theta_{\mu}$.
\end{remark}

\begin{remark}   \label{rem:6.7}
For $\mu \in \cG_1$ one has a natural way of extending the character
$\theta_{\mu} : \sym \to \bC$ to a unital homomorphism 
$\Theta_{\mu} : \sym \, [[z]] \to \bC \, [[z]]$, where a series 
$\uU (z) = \sum_{n=0}^{\infty} u_n z^n \in \sym \, [[z]]$ is mapped
by $\Theta_{\mu}$ to the series 
$\sum_{n=0}^{\infty} \theta_{\mu} (u_n) z^n \in \bC \, [[z]]$. 
It is useful to record what $\Theta_{\mu}$ does to the special 
series $\uE, \uH, \uP \in \sym \, [[z]]$ from Equations
(\ref{eqn:6.31})--(\ref{eqn:6.33}) of Remark \ref{rem:6.3}. First of 
all, from the very definition of $\Theta_{\mu}$ it is clear that
\begin{equation}   \label{eqn:6.6001}
\Theta_{\mu}  ( \uH ) = S_{\mu}.
\end{equation}
Then from Equation (\ref{eqn:6.34}) written in the form 
$\uE = 1 / \uH$ and the fact that $\Theta_{\mu}$ is a unital 
homomorphism, we infer that 
\begin{equation}  \label{eqn:6.6002}
\Theta_{\mu} ( \uE ) = 1/S_{\mu}.
\end{equation}
Finally, let us observe that when taking the logarithm of the 
left-hand side of (\ref{eqn:6.6002}), the log may be as well written
inside the brackets; this is because 
$\log \bigl( \, \Theta_{\mu} ( \uE ) \, \bigr)$ and 
$\Theta_{\mu} ( \log \uE )$ are both defined by using the 
Taylor series of the logarithm, and because of the homomorphic 
properties of $\Theta_{\mu}$. When 
we also invoke Equation (\ref{eqn:6.35}), we thus find that
\begin{equation}   \label{eqn:6.6003}
\Theta_{\mu}  ( \uP ) 
= \Theta_{\mu}  ( \log \uE ) 
= \log \bigl( \, \Theta_{\mu}  ( \uE ) \, \bigr) 
= \log (1/S_{\mu}) = - \log S_{\mu}.
\end{equation}
\end{remark}

\begin{definition}   \label{def:6.8}

$1^o$ For every $n \geq 2$ we denote
\begin{equation}   \label{eqn:6.81}
y_n = \sum_{ \begin{array}{c}
{\scriptstyle \pi = \{ A_1, \ldots , A_q \} } \\
{\scriptstyle \mathrm{in} \ NC(n-1) }
\end{array} } \ e_{|A_1|} \cdots e_{|A_q|} \in \sym .
\end{equation}
Thus for instance for $2 \leq n \leq 5$ we have
\begin{equation}   \label{eqn:6.82}
\left\{   \begin{array}{l}
y_2 = e_1, \  y_3 = e_2 + e_1^2, \
y_4 = e_3 + 3 e_1 e_2 + e_1^3, \\
                               \\
y_5 = e_4 + 4 e_1 e_3 + 2 e_2^2 + 6 e_1^2 e_2 + e_1^4.
\end{array}  \right.
\end{equation}
Clearly, every $y_n$ is a homogeneous symmetric function of degree
$n-1$.

$2^o$ We denote by $\Phi$ the unital algebra homomorphism from 
$\cYone ( \, = \bC [ Y_2, Y_3, \ldots \, ] \, )$ to $\sym$ which is 
uniquely determined by the requirement that $\Phi (Y_n) = y_n$,
$\forall \, n \geq 2$. 
\end{definition}

\begin{proposition} \label{prop:6.9}
For every $\mu \in \cG_1$ we have that
$\chi_{\mu} = \theta_{\mu}  \circ \Phi$, 
where $\chi_{\mu} \in \bX ( \cYone )$ and 
$\theta_{\mu} \in \bX ( \sym )$ are as in Definition
\ref{def:1.1} and Definition \ref{def:6.4}, respectively.
\end{proposition}

\begin{proof}
Both $\chi_{\mu}$ and $\theta_{\mu} \circ \Phi$ are characters of
$\cYone$, so in order to show that they are equal to each other it
suffices to verify that they agree on every $Y_n$, $n \geq 2$.

Let us consider the $R$-transform of $\mu$,
\[
R_{\mu} (z) = z + \sum_{n=2}^{\infty} \alpha_n z^n,
\]
and let us recall that $\chi_{\mu}$ is defined via
the requirement that $\chi_{\mu} (Y_n) = \alpha_n$, $\forall \, n
\geq 2$. Hence what we have to show is that $( \theta_{\mu} \circ
\Phi ) (Y_n) = \alpha_n$, or equivalently, that
\begin{equation}     \label{eqn:6.72}
\theta_{\mu} (y_n) = \alpha_n, \ \ \forall \, n \geq 2.
\end{equation}

On the other hand let us also consider the series
\[
1/S_{\mu} (z) =: 1 + \sum_{n=1}^{\infty} \gamma_n z^n.
\]
In Remark \ref{rem:6.7} it was observed that $\Theta_{\mu} ( \uE )
= 1/S_{\mu}$; when we write this equation in coefficients, it says
that $\theta_{\mu} (e_n) = \gamma_n$, $\forall \, n \geq 1$.
Consequently, when we apply the linear and multiplicative functional
$\theta_{\mu}$ to both sides of Equation (\ref{eqn:6.81}) we find
that
\[
\theta_{\mu} (y_n) = \sum_{ \begin{array}{c}
{\scriptstyle \pi = \{ A_1, \ldots , A_q \} } \\
{\scriptstyle \mathrm{in} \ NC(n-1) }
\end{array} } \ \gamma_{|A_1|} \cdots \gamma_{|A_q|},
\]
and (\ref{eqn:6.72}) follows from Lemma \ref{lemma:6.2}.
\end{proof}

Based on Proposition \ref{prop:6.9}, we can now place the
one-dimensional instance of the $\ls$-transform in the framework of
$\sym$.

\begin{proposition}    \label{prop:6.10}
For a distribution $\mu \in \cG_1$, the $\ls$-transform $\ls_{\mu}$
can be written in the form
\begin{equation}   \label{eqn:6.101}
\ls_{\mu} (z) = \sum_{n=2}^{\infty} \Bigl( \, ( \log \theta_{\mu} )
(y_n) \, \Bigr) z^n,
\end{equation}
where $\theta_{\mu} \in \bX ( \sym )$ is as in Definition
\ref{def:6.4}.
\end{proposition}

\begin{proof}
As pointed out in Equation (\ref{eqn:1.14}) of the introduction, in
the one-dimensional case the definition of the $\ls$-transform
simplifies to
\[
\ls_{\mu} (z) = \sum_{n=2}^{\infty} \Bigl( \, ( \log \chi_{\mu} )
(Y_n) \, \Bigr) z^n.
\]
Hence proving (\ref{eqn:6.101}) amounts to verifying that
$( \log \theta_{\mu} ) (y_n) = ( \log \chi_{\mu} ) (Y_n)$, 
$\forall \, n \geq 2$.
Clearly, this will follow if we can prove that 
$( \log \theta_{\mu} ) \circ \Phi = \log \chi_{\mu}$,
where $\Phi : \cYone \to \sym$ is as in Definition \ref{def:6.8}.2.
Since Proposition \ref{prop:6.9} gives us that 
$\theta_{\mu} \circ \Phi = \chi_{\mu}$, it will thus 
suffice to verify that for every character 
$\xi \in \bX ( \sym )$ we have
\begin{equation}   \label{eqn:6.104}
( \log \xi ) \circ \Phi = \log ( \xi \circ \Phi )
\end{equation}
(equality of linear functionals on $\cYone$).

In order to prove (\ref{eqn:6.104}), we start from the following
fact:
\begin{equation}   \label{eqn:6.105}
\left\{ \begin{array}{c}
\bX ( \sym ) \ni \xi \mapsto \xi \circ \Phi \in \bX ( \cYone ) \\
\mbox{ is a group homomorphism }
\end{array}  \right\}  .
\end{equation}
This holds because, as shown in Proposition \ref{prop:6.9}, the 
map in (\ref{eqn:6.105}) is the one which connects the two group 
isomorphisms $\mu \mapsto \chi_{\mu}$ and 
$\mu \mapsto \theta_{\mu}$ going from $\cG_1$ to $\bX ( \cYone )$ 
and to $\bX ( \sym )$, respectively. 
\footnote{We emphasize that at this point we are not using any 
properties that $\Phi$ might have in connection to the coalgebra 
structures of $\cYone$ and $\sym$; for the time being $\Phi$ is 
merely a unital algebra homomorphism.}
From (\ref{eqn:6.105}) we infer that for every 
$\xi \in \bX ( \sym )$ and every $\ell \geq 1$ we have
\begin{equation}   \label{eqn:6.86}
\bigl( \ee_{ { }_{\mathrm{Sym}} } - \xi \bigr)^{\ell} \circ \Phi 
= \bigl( \ee_{\cYone} - ( \xi \circ \Phi ) \bigr)^{\ell}
\end{equation}
(equality of functionals on $\cYone$, where 
$\ee_{ { }_{\mathrm{Sym}} }$ and 
$\ee_{\cYone}$ are counits, and the convolution powers are 
taken in the suitable algebras of functionals). Indeed, we have
\begin{align*}
( \ee_{ { }_{\mathrm{Sym}} } - \xi )^{\ell} \circ \Phi 
& = \Bigl( \sum_{j=0}^{\ell} { \ell \choose j } (-1)^j
\xi^j \Bigr) \circ \Phi                                    \\
& = \sum_{j=0}^{\ell} { \ell \choose j }(-1)^j 
( \xi^j \circ \Phi )                                       \\
& = \sum_{j=0}^{\ell} { \ell \choose j }(-1)^j 
( \xi \circ \Phi )^j \ \ \mbox{ (by (\ref{eqn:6.105})) }    \\
& = (\ee_{\cYone} - ( \xi \circ \Phi) \bigr)^{\ell} .
\end{align*}
Finally we are just left to observe that for every 
$\xi \in \bX ( \sym )$ we have
\begin{align*}
( \log \xi ) \circ \Phi
& = \Bigl( - \sum_{\ell =1}^{\infty} \frac{1}{\ell}
( \ee_{ { }_{\mathrm{Sym}} } - \xi )^{\ell} \Bigr) \circ\Phi   \\
& = - \sum_{\ell =1}^{\infty} \frac{1}{\ell} 
\bigl( \ee_{\cYone} -\ ( \xi \circ \Phi ) \Bigr)^{\ell}
\ \ \mbox{ (by (\ref{eqn:6.86})) }                           \\
& = \log ( \xi \circ \Phi ),
\end{align*}
as claimed in (\ref{eqn:6.104}).
\end{proof}

It is now easy to obtain the precise formula relating the
$\ls$-transform with the log of the $S$-transform -- we will only
have to combine Proposition \ref{prop:6.10} with some basic
properties of the infinitesimal characters of $\sym$, which are
recorded in the next lemma.

\begin{lemma}   \label{lemma:6.11}

$1^o$ Let $\xi$ be an infinitesimal character in $\bI ( \sym )$, and
let $u_1 , u_2$ be homogeneous symmetric functions of degrees $\geq
1$. Then $\xi (u_1 u_2) = 0$.

$2^o$ Let $\eta$ be a character in $\bX ( \sym )$, and consider its
logarithm $\xi := \log \eta \in \bI ( \sym )$. Then 
\begin{equation}    \label{eqn:6.91}
\xi ( p_n ) = \eta ( p_n ), \ \ \forall \, n \geq 1.
\end{equation}

$3^o$ Let $\xi$ be an infinitesimal character in $\bI ( \sym )$. Then
we have that
\begin{equation}    \label{eqn:6.92}
\xi ( y_n ) = \frac{(-1)^n}{n-1} \xi ( p_{n-1} ), \ \ \forall 
\, n \geq 2.
\end{equation}
\end{lemma}

\begin{proof}

$1^o$ We have that 
\begin{align*}
\xi (u_1 u_2)
& = \xi (u_1) \ee (u_2) + \ee (u_1) \xi (u_2)
\ \ \mbox{ (by the definition of $\bI ( \sym )$) }            \\
& = \xi (u_1) \cdot 0 + 0 \cdot \xi (u_2) 
\ \ \mbox{ (by the definition of $\ee$) }                     \\
& =0.
\end{align*}

$2^o$ Note that $( \ee -  \eta )( 1 )= 1-1=0$, which implies 
that for every $\ell \geq 1$ we have
\begin{eqnarray*}
(\ee -  \eta )^{\ell} (1) = 
( (\ee - \eta ) \otimes \cdots \otimes ( \ee - \eta ) )
(1 \otimes \cdots \otimes 1) = 0.
\end{eqnarray*}
For $\ell \geq 2$ this implies that
\[
(\varepsilon - \eta )^{\ell} (p_n) 
= ((\varepsilon -  \eta )^{\ell -1}\otimes (\varepsilon - \eta ))
(p_n\otimes 1+1\otimes p_n) = 0
\]
(where we wrote $( \varepsilon -  \eta )^{\ell} =
( \varepsilon -  \eta )^{\ell - 1} \cdot ( \varepsilon - \eta )$, 
and we used formula (\ref{eqn:6.38}) for $\Delta (p_n)$). Hence 
indeed, we find that 
\[
\xi (p_n) = - \sum_{\ell =1}^{\infty} \frac{1}{\ell} 
\Bigl( \, ( \ee - \eta )^{\ell} \, \Bigr) (p_n) 
= - ( \ee - \eta ) (p_n) = \eta (p_n).
\]

$3^o$ Equation (\ref{eqn:6.81}) which defines $y_n$ can also be 
written in the form
\[
y_n = e_{n-1} + \sum_{ \begin{array}{c}
{\scriptstyle \pi = \{ A_1, \ldots , A_q \} } \\
{\scriptstyle \mathrm{in} \ NC(n-1), \ q \geq 2 }
\end{array} } \ e_{|A_1|} \cdots e_{|A_q|} ;
\]
hence by using part $1^o$ of the lemma we find that
\begin{equation}   \label{eqn:6.94}
\xi (y_n) = \xi ( e_{n-1} ), \ \ \forall \, n \geq 2.
\end{equation}
Now let us look again at Equation (\ref{eqn:6.35}) which connects 
the series $\uE, \uP \in \sym \, [[z]]$. By taking derivative with 
respect to $z$ we obtain that $\uP ' \cdot \uE = \uE '$, and when 
written coefficientwise this says that
\begin{equation}   \label{eqn:6.95}
p_1 e_{n-1} - p_2 e_{n-2} + \cdots + (-1)^n p_{n-1} e_1 + (-1)^{n+1}
p_n = n e_n, \ \ \forall \, n \geq 1.
\end{equation}
We apply $\xi$ to both sides of (\ref{eqn:6.95}) and
invoke part $1^o$ of the lemma to conclude that 
\begin{equation}   \label{eqn:6.96}
(-1)^{n+1} \xi (p_n) = n \xi ( e_n ), \ \ \forall \, n \geq 1.
\end{equation}
Clearly, the required Equation (\ref{eqn:6.92}) follows from
(\ref{eqn:6.94}) and (\ref{eqn:6.96}).
\end{proof}

\begin{corollary}     \label{cor:6.12}
Let $\mu$ be a distribution in $\cG_1$, and consider the power
series $S_{\mu}, \ls_{\mu} \in \bC \, [[z]]$. Then 
$\ls_{\mu} (z) = -z \log S_{\mu} (z)$.
\end{corollary}

\begin{proof} The constant and the linear term vanish in both
$\ls_{\mu} (z)$ and $-z \log S_{\mu} (z)$, hence we only have to 
look at coefficients of $z^n$ for $n \geq 2$.

So let us fix an integer $n \geq 2$. By Proposition \ref{prop:6.10}, 
the coefficient of $z^n$ in $\ls_{\mu}$ is equal to 
$( \log \theta_{\mu}) (y_n)$. But on the other hand we have that
\begin{align*}
( \log \theta_{\mu} ) (y_n) 
& = \frac{(-1)^n}{n-1} ( \log \theta_{\mu}) ( p_{n-1} )
\mbox{  (by Lemma \ref{lemma:6.11}.3) }                            \\
& = \frac{(-1)^n}{n-1} \theta_{\mu} ( p_{n-1} )
\mbox{  (by Lemma \ref{lemma:6.11}.2) }                            \\
& =   \theta_{\mu} \mbox{ (coefficient of $z^{n-1}$ in $\uP$)}   
\mbox{ (by (\ref{eqn:6.33})) }                                    \\
& =   \mbox{coefficient of $z^{n-1}$ in $\Theta_{\mu} ( \uP )$}   \\
& =   \mbox{coefficient of $z^{n-1}$ in $( -\log S_{\mu} )$, }
\end{align*}
where at the last two equality signs we used the homomorphism
$\Theta_{\mu}$ from Remark \ref{rem:6.7}, and its property
observed in Equation (\ref{eqn:6.6003}). So we obtained equality
with the coefficient of $z^n$ in $-z \log S_{\mu} (z)$, as required.
\end{proof}

\vspace{2.5cm}

\begin{center}
{\bf\large 7. Some remarks about \boldmath{$\cYone$} }
\end{center}
\setcounter{section}{7}
\setcounter{equation}{0}
\setcounter{theorem}{0}

The goal of the present section is to clarify a few points concerning 
the Hopf algebra $\cYone$, that arise naturally from the 
considerations in Section 6. 

We will continue to use the same framework and notations as in 
Section 6. Recall that the connection between $\cYone$ and 
symmetric functions was formally achieved by the homomorphism 
$\Phi : \cYone \to \sym$ from Definition \ref{def:6.8}.2. 
In Section 6 it came in handy to introduce and use $\Phi$ 
only in its capacity of homomorphism of unital algebras; 
but $\Phi$ is in fact more than that -- it is an isomorphism of graded 
Hopf algebras. This will be proved in Theorem \ref{thm:7.5} below.
Our method of proof will be to observe that $\cYone$ is cocommutative, 
and then to take advantage of the special properties that primitive 
elements are known to have in cocommutative Hopf algebras.

\begin{lemma}  \label{lemma:7.1}
The Hopf algebra $\cYone$ is cocommutative.
\end{lemma}

\begin{proof} We have to show that $F \circ \Delta = \Delta$, where 
$F : \cYone \otimes \cYone \to  \cYone \otimes \cYone$ is the flip 
automorphism acting by $F( P \otimes Q) = Q \otimes P$ for
$P,Q \in \cYone$. Since both $\Delta$ and $F \circ \Delta$ are unital
algebra homomorphisms from $\cYone$ to $\cYone \otimes \cYone$, it 
suffices to verify that they agree on the set of generators
$\{ Y_n \mid n \geq 2 \}$ of $\cYone$.

So let us fix an $n \geq 2$. The formula defining $\Delta$ (see
Equation (\ref{eqn:3.21}) in Definition \ref{def:3.2}.1) gives us here
that
\begin{equation}  \label{eqn:7.11}
\Delta ( Y_n ) = \sum_{\pi \in NC(n)}
Y_{\pi} \otimes Y_{ K( \pi ) },
\end{equation}
where for $\pi = \{ A_1, \ldots , A_q \} \in NC(n)$ we put
\begin{equation}  \label{eqn:7.12}
Y_{\pi} := Y_{ |A_1| } \cdots Y_{ |A_q| }.
\end{equation}
On the other hand we have
\begin{align*}
( F \circ \Delta ) (Y_n)
& = F \bigl( \, \sum_{\pi \in NC(n)}
Y_{\pi} \otimes Y_{ K( \pi ) } \, \bigr)                     \\
& = \sum_{\pi \in NC(n)} Y_{ K( \pi ) } \otimes Y_{ \pi }    \\ 
& = \sum_{\rho \in NC(n)} Y_{ \rho } \otimes Y_{ K^{-1} ( \rho ) }   
\ \ \mbox{ (by substitution $K( \pi ) = \rho$). }
\end{align*}
Now, it is immediately seen that for every $\rho \in NC(n)$, the
partitions $K( \rho )$ and $K^{-1} ( \rho )$ have the same block 
structure -- that is, they can be written as
\begin{equation}   \label{eqn:7.13}
\left\{  \begin{array}{c}
K( \rho ) = \{ B_1, \ldots , B_q \}, \ \ 
K^{-1} ( \rho ) = \{ C_1, \ldots , C_q \},               \\
\mbox{ with } |B_1| = |C_1|, \ldots , |B_q| = |C_q|.
\end{array}  \right.
\end{equation}
(See Exercise 9.23 on page 147 of \cite{NS06}; or directly from 
formula (\ref{eqn:2.11}) in Notation \ref{def:2.1}, observe that 
the associated permutations
$P_{K( \rho )} = P_{\rho}^{-1} P_{1_n}$ and
$P_{K^{-1} ( \rho )} = P_{1_n} P_{\rho}^{-1}$ belong to the same
conjugacy class of the symmetric group.) From (\ref{eqn:7.13}) it is
immediate that $Y_{ K^{-1} ( \rho )} = Y_{K ( \rho )}$, 
$\forall \, \rho \in NC(n)$, so we find that 
\[
( F \circ \Delta ) (Y_n)  
= \sum_{\rho \in NC(n)} Y_{ \rho } \otimes Y_{ K ( \rho ) }   
= \Delta (Y_n),
\]
as required.
\end{proof}

\begin{remark}  \label{rem:7.2}
Consider the set of primitive elements of $\cYone$,
\begin{equation}   \label{eqn:7.21}
\PrimY := \{ P \in \cYone \mid \Delta (P) = P \otimes 1 
+ 1 \otimes P \}.
\end{equation}
The commutator operation $[P,Q] := PQ-QP$ turns $\PrimY$ into
a Lie algebra, and a fundamental result from the theory of 
cocommutative Hopf algebras (the Cartier-Kostant-Milnor-Moore
theorem, see e.g. Theorem 13.0.1 on p. 274 of \cite{S69})
allows us to identify $\cYone$ as the universal envelopping 
algebra of $\PrimY$. For our discussion here, geared towards the
proof of Theorem \ref{thm:7.5} below, it will suffice to use the 
weaker statement that $\PrimY$ generates $\cYone$ as a unital 
algebra -- see Corollary 13.0.3 on p. 278 of \cite{S69}. Since 
it is immediate that every primitive element can be written as 
a sum of homogeneous primitive elements, the statement about 
primitives that we want to retain can thus be recorded as follows:
\begin{equation}   \label{eqn:7.22}
\left\{  \begin{array}{c}
\mbox{ $\cYone$ is generated as a unital algebra}   \\
\mbox{ by the set of its homogeneous primitive elements }
\end{array}  \right\}  .
\end{equation}

Before proceeding to the Theorem \ref{thm:7.5} we will 
also record, in Lemma \ref{lemma:7.4}, a simple observation about
homogeneous symmetric functions.
\end{remark}

\begin{notation}  \label{def:7.3}
Let $n$ be a positive integer.

$1^o$ We use the standard notation ``$\lambda \vdash n$'' to denote a
partition $\lambda$ of the number $n$ (that is, $\lambda$ is a tuple
$(i_1, i_2, \ldots , i_m)$ of integers 
$i_1 \geq i_2 \geq \cdots \geq i_m >0$ with 
$i_1 + i_2 + \cdots + i_m =n$). For 
$\lambda = (i_1, \ldots , i_m) \vdash n$, the number $m$ of parts of 
$\lambda$ is called the {\em length} of $\lambda$,
and is denoted by $\ell ( \lambda )$.

$2^o$ We denote by $\sym_{\, n}$ the vector space of homogeneous 
symmetric functions of degree $n$. In $\sym_{\, n}$ we will use the 
linear basis $\{ p_{\lambda} \mid \lambda \vdash n \}$, where for 
every $\lambda = ( i_1, i_2, \ldots , i_m ) \vdash n$ we put
\begin{equation}  \label{eqn:7.31}
p_{\lambda} := p_{i_1} p_{i_2} \cdots p_{i_m} \in \sym_{\, n}.
\end{equation}
(For the verification that $\{ p_{\lambda} \mid \lambda \vdash n \}$
is indeed a basis for $\sym_{\, n}$, see e.g. Corollary 7.7.2 on 
p. 298 of \cite{S99}.)
\end{notation}

\begin{lemma}  \label{lemma:7.4}
Let $n$ be a positive integer and let $u$ be a symmetric function 
in $\sym_{\, n}$ which has the property that
\begin{equation}   \label{eqn:7.41}
\xi^2 (u) = 2 \xi (u), \ \ \forall \, \xi \in \bX ( \sym ).
\end{equation}
Then $u \in \bC p_n$.
\end{lemma}

\begin{proof} Consider the writing of $u$ in terms of the basis 
$\{ p_{\lambda} \mid \lambda \vdash n \}$:
\begin{equation}   \label{eqn:7.42}
u = \sum_{\lambda \vdash n} \, t_{\lambda} p_{\lambda}, \ \ 
\mbox{ with $( t_{\lambda} )_{\lambda \vdash n}$ from $\bC$.}
\end{equation}
We leave it as an exercise to the reader to use the formula
defining $\xi^2$ and the fact that each of $p_1, p_2, \ldots , p_n$
is a primitive element of $\sym$ in order to verify that 
(\ref{eqn:7.42}) implies:
\begin{equation}   \label{eqn:7.43}
\xi^2 (u) = 2 \xi (u)  
+ \sum_{ \begin{array}{c}
{\scriptstyle \lambda \vdash n,} \\
{\scriptstyle \ell ( \lambda ) \geq 2}
\end{array}  } \  t_{\lambda} N( \lambda ) \xi ( p_{\lambda} ), 
\ \ \forall \, \xi \in \bX ( \sym ),
\end{equation}
where for every $\lambda \vdash n$ with $\ell ( \lambda ) \geq 2$,
$N ( \lambda )$ is a positive integer depending only on $\lambda$.
(The combinatorial significance of $N( \lambda )$ is that it counts
in how many ways $\lambda$ can be obtained by merging together two
partitions $\lambda ' \vdash  m$ and $\lambda '' \vdash n-m$ for some
$1 \leq m < n$.) The hypothesis given on $u$ thus amounts to the 
fact that 
\begin{equation}   \label{eqn:7.44}
\sum_{ \begin{array}{c}
{\scriptstyle \lambda \vdash n,} \\
{\scriptstyle \ell ( \lambda ) \geq 2}
\end{array}  } \  t_{\lambda} N( \lambda ) \xi ( p_{\lambda} )
= 0, \ \ \forall \, \xi \in \bX ( \sym ).
\end{equation}
But a character $\xi \in \bX ( \sym )$ may take any prescribed set 
of values on $p_1, p_2, \ldots , p_n$; by using this fact it is easy
to see that (\ref{eqn:7.44}) can only hold if 
$t_{\lambda} N( \lambda ) =0$ (hence $t_{\lambda} = 0$) for every 
$\lambda \vdash n$ with $\ell ( \lambda ) \geq 2$, and
the conclusion that $u \in \bC p_n$ follows.
\end{proof}

\begin{theorem}   \label{thm:7.5}
The map $\Phi : \cYone \to \sym$ from Definition \ref{def:6.8}.2 is 
an isomorphism of graded Hopf algebras.
\end{theorem}

\begin{proof} $\Phi$ is by definition a unital algebra homomorphism.
Directly from the definition it is also immediate that $\Phi$ 
respects the gradings considered on $\cYone$ and on $\sym$, and that
$\Phi$ respects counits (we have
$\ee_{ { }_{\mathrm{Sym}} } \circ \Phi = \ee_{\cYone}$). 

We next verify that $\Phi$ is surjective. It clearly suffices to
verify that $\Phi ( \cYone ) \ni e_n$ for every $n \geq 1$, which 
we do by induction on $n$. The base case $n=1$ holds because 
$e_1 = y_2 = \Phi ( Y_2 )$. For the induction step let us assume that
$e_1, \ldots , e_{n-1}$ belong to $\Phi ( \cYone )$, and let us prove 
that $e_n$ belongs to $\Phi ( \cYone )$ as well, for some $n \geq 2$. 
To this end we observe that the formula used to define $y_{n+1}$ in 
Definition \ref{def:6.8}.1 can be written in the form
\begin{equation}   \label{eqn:7.51}
e_n = y_{n+1} - \sum_{ \begin{array}{c}
{\scriptstyle \pi = \{ A_1, \ldots , A_q \} \ \mathrm{in} } \\
{\scriptstyle NC(n), \ \mathrm{with} \ q \geq 2}
\end{array} } \ e_{|A_1|} \cdots e_{|A_q|}.
\end{equation}
The right-hand side of (\ref{eqn:7.51}) is a polynomial expression
in $y_{n+1}$ and $e_1, \ldots , e_{n-1}$, which all belong to
$\Phi ( \cYone )$. Since $\Phi ( \cYone )$ is a subalgebra of 
$\sym$, we thus conclude that $e_n \in \Phi ( \cYone )$, as 
required.

Now, in the context where we know that $\Phi$ is surjective and 
respects the gradings, the injectivity of $\Phi$ will also follow 
if we can check that for every $n \geq 0$ the homogeneous spaces 
$\cYone_n \subseteq \cYone$ and $\sym_{\, n} \subseteq \sym$ have
the same finite dimension. And indeed,
$\dim ( \cYone_0 ) = \dim ( \sym_{ \, 0} ) = 1$, while for every 
$n \geq 1$ we have 
\begin{equation}   \label{eqn:7.52}  
\dim ( \sym_{ \, n} ) = 
\ \vert \, \{ \lambda \mid \lambda \vdash n \} \, \vert \
= \dim ( \cYone_n ).
\end{equation}
The first of the two equalities in (\ref{eqn:7.52}) follows from the
fact that $\{ p_{\lambda} \mid \lambda \vdash n \}$ is a basis for
$\sym_{\, n}$, while the second is an immediate consequence of
the fact that $\cYone$ is just $\bC [Y_2, Y_3, \ldots \, ]$, with
degree$(Y_m) = m-1$ for every $m \geq 2$. This completes the 
verifications that $\Phi$ is bijective.

We are left to prove that $\Phi$ respects the comultiplications of
$\cYone$ and $\sym$, i.e. that 
\begin{equation}   \label{eqn:7.53}  
\Delta_{ { }_{\mathrm{Sym}} } \circ \Phi =
( \Phi \otimes \Phi ) \circ \Delta_{\cYone}.
\end{equation}
Due to (\ref{eqn:7.22}) of Remark \ref{rem:7.2} and the fact that
both sides of Equation (\ref{eqn:7.53}) are unital homomorphisms 
from $\cYone$ to $\sym \otimes \sym$, it is sufficient to verify that
\[
\Bigl( \, \Delta_{ { }_{\mathrm{Sym}} } \circ \Phi \, \Bigr) (P) =
\Bigl( \, ( \Phi \otimes \Phi ) \circ \Delta_{\cYone} \, \Bigr) (P)
\]
when $P$ is a homogeneous primitive element of $\cYone$.
So let $P$ be a such an element, of degree $n$.
We claim that the symmetric function 
$u := \Phi (P) \in \sym_{\, n}$ satisfies the hypothesis of 
Lemma \ref{lemma:7.4}. Indeed, for every character 
$\xi \in \bX ( \sym )$ we write:
\begin{align*} 
\xi^2 (u) 
& = ( \xi^2 \circ \Phi ) (P)                                    \\
& = ( \xi \circ \Phi )^2 (P)  \ \ \mbox{ (by (\ref{eqn:6.105}) 
                      in proof of Proposition \ref{prop:6.10})} \\
& = \Bigl( ( \xi \circ \Phi ) \otimes ( \xi \circ \Phi ) \Bigr)
    (P \otimes 1 + 1 \otimes P)                                 \\
& = 2 \xi (u).
\end{align*}
Hence Lemma \ref{lemma:7.4} applies and gives us that $u$ is 
a scalar multiple of the power sum symmetric function $p_n$. This 
implies in particular that $u$ is primitive in $\sym$, and it 
follows that
\[
\Bigl( \, \Delta_{ { }_{\mathrm{Sym}} } \circ \Phi \, \Bigr) (P) =
u \otimes 1 + 1 \otimes u =
\Bigl( \, ( \Phi \otimes \Phi ) \circ \Delta_{\cYone} \, \Bigr) (P),
\]
as we wanted.
\end{proof}

\begin{remark}  \label{rem:7.6}
In the proof of Theorem \ref{thm:7.5}, an alternative approach to
the fact that $\Phi$ respects comultiplication would go by proving, 
directly from the definitions, that every symmetric function $y_n$
comultiplies in $\sym$ by the same formula as the one used for the
comultiplication of $Y_n \in \cYone$. That is, one could go for a 
direct proof of the formula:
\begin{equation}    \label{eqn:7.61}
\Delta (y_n) = \sum_{\pi \in NC(n)} \ y_{\pi} \otimes y_{K( \pi )},
\ \ \forall \, n \geq 2,
\end{equation}
where for $\pi = \{ A_1, A_2, \ldots , A_q \} \in NC(n)$ we put
\begin{equation}    \label{eqn:7.62}
y_{\pi} := y_{ |A_1| } \  y_{ |A_2| } \cdots y_{ |A_q| }
\in \sym ,
\end{equation}
with the convention that $y_1 := 1$. By ``direct proof'' we mean 
here a direct derivation of (\ref{eqn:7.61}) from 
Equation (\ref{eqn:6.81}) which defines the $y_n$'s in terms of 
$e_n$'s, combined with the simple formula (\ref{eqn:6.37}) for the
comultiplication of $e_n$. (It is instructive to see how this works
for some concrete small values of $n$, e.g. pick $n=3$ and verify 
directly that $\Delta ( y_3 )$ =
$y_3 \otimes 1 + 3 y_2 \otimes y_2 + 1 \otimes y_3$.)

For the readers who want to try their hand at the direct derivation
of (\ref{eqn:7.61}), we offer the following tip: the formula 
defining $y_n$ in (\ref{eqn:6.81}) can be re-written in the form of
a recursion,
\begin{equation}    \label{eqn:7.63}
y_n = \sum_{m=2}^n  \Bigl( \, e_{m-1} \cdot
\sum_{1= i_1<i_2< \cdots < i_m =n} \
y_{i_2 - i_1} \, y_{i_3 - i_2} \cdots y_{i_m - i_{m-1}} \, \Bigr) .
\end{equation}
This recursion is essentially equivalent to the ``functional equation
of the $R$-transform'' from Lecture 16 of \cite{NS06}, which was also
invoked in the proof of Lemma \ref{lemma:6.2}. One can use 
(\ref{eqn:7.63}) and the formula (\ref{eqn:6.37}) for the
comultiplication of the $e_n$'s in order to prove {\em by induction}
on $n$ that (\ref{eqn:7.61}) holds. Proceeding by induction somewhat
simplifies the argument, but this ``direct proof'' that $\Phi$ 
respects comultiplication still remains (by quite a bit) more involved
than the approach which we chose to present in the proof shown above 
for Theorem \ref{thm:7.5}.
\end{remark}

\vspace{10pt}

Another way of looking at the homomorphism $\Phi : \cYone \to \sym$
is by placing it in the framework of combinatorial Hopf algebras 
from \cite{ABS06}. We will explain how this goes in the next remark,
where it is convenient to start by pointing out a connection which 
$\cYone$ has with the theory of incidence Hopf algebras.

\begin{remark} \label{rem:7.7}

$1^o$ Let $\cR$ be Rota's Hopf algebra of isomorphism classes of
finite graded posets. We use the name ``Rota's Hopf algebra'' and 
the notation $\cR$ by following Example 2.2 of \cite{ABS06}. Standard
references for this graded connected Hopf algebra are \cite{JR79}, 
\cite{S94}. Here is a brief summary of facts about $\cR$ that we
want to refer to. In (i)-(iv) below, by ``finite graded poset'' we 
mean a finite poset $G$ with minimum element $0_G$ and maximum 
element $1_G$, such that all the saturated chains in $G$ have the 
same length $r$; the number $r$ is called the rank of the poset.

(i) For every finite graded poset $G$ one has a homogeneous 
element $\overline{G} \in \cR$, of degree equal to the rank of $G$.

(ii) If $G_1, \ldots , G_n$ are finite graded posets such that 
$G_i \not\simeq G_j$ for $1 \leq i<j \leq n$, then the elements 
$\overline{G_1}, \ldots , \overline{G_n}$ of $\cR$ are linearly 
independent.

(iii) The multiplication of $\cR$ is such that
$\overline{G_1} \cdot \overline{G_2} = \overline{ G_1 \times G_2 }$,
where $G_1 \times G_2$ denotes the direct product of the finite 
graded posets $G_1$ and $G_2$.

(iv) The comultiplication of $\cR$ is such that for every finite 
graded poset $G$ one has
\[
\Delta ( \, \overline{G} \, ) = \sum_{x \in G}
\overline{ [ 0_G, x ] } \otimes \overline{ [ x, 1_G ] },
\]
where we denote
$[ 0_G, x ] := \{ y \in G \mid y \leq x \}$ and
$[ x, 1_G ] := \{ y \in G \mid y \geq x \}$, $x \in G$.

\noindent
For a more detailed presentation of how $\cR$ is constructed, see 
e.g. Section 3 of \cite{S94}.

What makes $\cR$ come into our considerations is that one has a 
natural embedding of $\cYone$ into $\cR$. More precisely, let 
$\Lambda : \cYone \to \cR$ be the unital algebra homomorphism 
uniquely determined by the requirement that
\begin{equation} \label{eqn:7.71} 
\Lambda (Y_n) = \overline{ NC(n) }, \ \ \forall \, n \geq 2.
\end{equation}
The assignment for $\Lambda (Y_n)$ in (\ref{eqn:7.71}) makes 
sense, as $NC(n)$ is a finite graded poset of rank $n-1$ (see e.g.
Exercise 10.29 on p. 171 of \cite{NS06}). We have that
$\Lambda$ is an injective homomorphism of graded Hopf algebras.
The verification of this statement is left as an exercise to the 
reader, we only indicate references for the two details of the 
verification that are non-trivial.

(a) Why is $\Lambda$ injective: in view of (ii) and (iii) above, this 
comes down to the fact that if 
$2 \leq i_1 \leq i_2 \leq \cdots \leq i_m$
and $2 \leq j_1 \leq j_2 \leq \cdots \leq j_n$  are such that 
\[
NC( i_1) \times \cdots \times NC(i_m) \simeq
NC( j_1) \times \cdots \times NC(j_n),
\]
then $m=n$ and $i_1 = j_1, \ldots , i_n = j_n$. For a proof of this 
fact, see Proposition 9.38 on p. 152 of \cite{NS06}.

(b) Why does $\Lambda$ respect comultiplication: in view of (iii) 
and (iv) above, this follows immediately from the fact that for
$\pi = \{ A_1, \ldots , A_q \} \in NC(n)$ with 
$K( \pi ) =: \{ B_1, \ldots , B_r \}$, the intervals 
$[0_n, \pi ]$ and $[ \pi, 1_n ]$ of $NC(n)$ satisfy
\[
\begin{array}{l}
{ [ 0_n , \pi ] } \simeq 
NC( \, |A_1| \, ) \times \cdots \times NC( \, |A_q| \, ),       \\
                                                                \\
{ [ \pi, 1_n ] } \simeq  { [ 0_n , K( \pi ) ] } \simeq
NC( \, |B_1| \, ) \times \cdots \times NC( \, |B_r| \, );
\end{array}
\]
see e.g. the discussion on pp. 149-150 of \cite{NS06}.

$2^o$ A benefit coming from the embedding of $\cYone$ into $\cR$ is 
that one gets a natural choice for a ``zeta character'' on $\cYone$. 
Indeed, $\cR$ has a zeta character $\zeta_{\cR}$ which acts by
\[
\zeta_{\cR} ( \, \overline{G} \, ) = 1, \ \ 
\mbox{ for every finite graded poset $G$; }
\]
see Example 2.2 of \cite{ABS06}. By restricting to $\cYone$ we thus 
see that the zeta character of $\cYone$ should act by the 
prescription that
\begin{equation}   \label{eqn:7.72}
\zeta_{\cYone} ( Y_n ) = 1, \ \ \forall \, n \geq 2.
\end{equation}

Now, a combinatorial Hopf algebra is defined in the paper 
\cite{ABS06} as a pair $( \cB , \zeta )$ where $\cB$ is a graded 
connected Hopf algebra and $\zeta$ is a character of $\cB$. So 
in particular $( \cYone , \zeta_{\cYone} )$ is 
a combinatorial Hopf algebra, where $\cYone$ is cocommutative. 
Theorem 4.3 of \cite{ABS06} asserts that if $( \cB , \zeta )$ is a 
combinatorial Hopf algebra such that $\cB$ is cocommutative, then
there exists a unique homomorphism of graded Hopf algebras 
$\Psi : \cB \to \sym$ such that 
$\zeta_{ { }_{\mathrm{Sym}} } \circ \Psi = \zeta$, where 
$\zeta_{ { }_{\mathrm{Sym}} }  \in \bX ( \sym )$ is the special 
character determined by the requirement that 
\begin{equation}   \label{eqn:7.73}
\zeta_{ { }_{\mathrm{Sym}} } ( e_1 ) = 1 \mbox{ and }
\zeta_{ { }_{\mathrm{Sym}} } ( e_n ) = 0 
\mbox{ for every $n \geq 2$. }
\end{equation}
We claim that in the particular case when $( \cB , \zeta )$ is 
$( \cYone , \zeta_{\cYone} )$, this gives us precisely the 
homomorphism $\Phi$ from Theorem \ref{thm:7.5}. That is, we have the 
following proposition.
\end{remark}

\begin{proposition}   \label{prop:7.8}
$\zeta_{ { }_{\mathrm{Sym}} } \circ \Phi = \zeta_{\cYone}$.
\end{proposition}

\begin{proof} For every $n \geq 2$ we have
\[
\Bigl( \, \zeta_{ { }_{\mathrm{Sym}} } \circ \Phi \, \Bigr) (Y_n)
= \zeta_{ { }_{\mathrm{Sym}} } (y_n) \ \ 
\mbox{ (since $\Phi ( Y_n ) = y_n$, by Definition \ref{def:6.8}.2) } 
\]
\begin{equation}   \label{eqn:7.81} 
= \sum_{  \begin{array}{c} 
{\scriptstyle \pi \in NC(n-1),}  \\
{\scriptstyle \pi = \{ A_1, \ldots , A_q \} }
\end{array} } \ \zeta_{ { }_{\mathrm{Sym}} } (e_{ |A_1| }) \cdots
\zeta_{ { }_{\mathrm{Sym}} } (e_{ |A_q| }) 
\ \ \mbox{ (by def. of $y_n$ in (\ref{eqn:6.81})).} 
\end{equation}
But from (\ref{eqn:7.73}) it is clear 
that the only non-zero term in the sum (\ref{eqn:7.81}) is the one 
corresponding to the partition $0_{n-1} \in NC(n-1)$, and this 
term is equal to 1. So we obtain that 
$\Bigl( \, \zeta_{ { }_{\mathrm{Sym}} } \circ \Phi \, \Bigr)$
$(Y_n) = 1 = \zeta_{\cYone} (Y_n)$, $\forall \, n \geq 2$, and the 
proposition follows.
\end{proof}

\begin{remark}   \label{rem:7.9}
Having placed the homomorphism $\Phi$ in the framework of 
combinatorial Hopf algebras leads to an interesting alternative 
description of the symmetric functions $\{ y_n \mid n \geq 2 \}$, 
as linear combinations of monomial quasi-symmetric functions. For 
every $m$-tuple of positive integers $(r_1, \ldots , r_m )$, the 
corresponding monomial quasi-symmetric function 
$M_{ (r_1, \ldots , r_m) }$ is defined as
\begin{equation}  \label{eqn:7.91}
M_{ (r_1, \ldots , r_m) } = 
\sum_{1 \leq i_1 < i_2 < \cdots < i_m} \
x_{i_1}^{r_1} x_{i_2}^{r_2} \cdots x_{i_m}^{r_m} ,
\end{equation}
where $\{ x_i \mid i \geq 1 \}$ is a family of commuting 
indeterminates. The elementary symmetric functions $e_n$ used
earlier in the paper are of course a particular case of 
(\ref{eqn:7.91}),
\begin{equation}  \label{eqn:7.92}
e_n = M_{ ( \ \underbrace{1, \ldots , 1}_{n} \ ) }, 
\ \ \forall \, n \geq 1.
\end{equation}
Hence one can take the point of view that Equation (\ref{eqn:6.81})
defines $y_n$ as a sum of products of monomial quasi-symmetric 
functions, and one can expand these products in order to obtain 
$y_n$ as a linear combination of $M_{ (r_1, \ldots , r_m) }$'s 
with $r_1 + \cdots + r_m = n-1$. To give a concrete example, for 
$n=4$ one has
\begin{align*}
y_4  
& = e_3 + 3 e_1 e_2 + e_1^3 \ \ \mbox{ (as in (\ref{eqn:6.82})) }   \\
& = M_{ (1,1,1) } + 3 M_{ (1) } M_{ (1,1) } + M_{ (1) }^3           \\
& = M_{ (1,1,1) } 
  + 3 \Bigl( \, 3M_{ (1,1,1) } + M_{ (2,1) } + M_{ (1,2) } \, \Bigr) 
  + \Bigl( \, 6M_{ (1,1,1) } + 3M_{ (2,1) } + 3M_{ (1,2) } 
                             + M_{ (3) } \, \Bigr)                  \\
& = 16M_{ (1,1,1) } + 6M_{ (2,1) } + 6M_{ (1,2) } + M_{ (3) }.
\end{align*}

In the above calculation for $n=4$ we ended up with a sum of 29
terms,which turn out to correspond to the $1+12+16$ chains in 
$NC(4)$ that were also used for illustration in Example \ref{ex:4.6}.
Indeed, the results about combinatorial Hopf algebras from
Theorems 4.1 and 4.3 of \cite{ABS06} contain a concrete recipe
for how to write $\Phi (Y_n)$ as a linear combination of monomial
quasi-symmetric functions. Due to the embedding of $\cYone$ into 
Rota's Hopf algebra from Remark \ref{rem:7.7}.1 this recipe reduces
to one previously found in Ehrenborg's paper \cite{E96} (see the 
discussion in Example 4.4 of \cite{ABS06}), and we get that
\begin{equation}    \label{eqn:7.93}
y_n = \sum_{ \begin{array}{c}
{\scriptstyle \Gamma = ( \pi_o, \pi_1, \ldots , \pi_{\ell} ) }     \\
{\scriptstyle \mathrm{chain} \ \mathrm{in} \ NC(n)}
\end{array}  } \ 
M_{ ( | \pi_0 | - | \pi_1 |,
| \pi_1 | - | \pi_2 |, \ldots ,
| \pi_{\ell -1} | - | \pi_{\ell} | ) }, \ \ \forall 
\, n \geq 2.
\end{equation}
It is likely that the equivalence between formulas (\ref{eqn:6.81}) 
and (\ref{eqn:7.93}) can be proved directly from the rules
of how products of monomial quasi-symmetric functions are formed, 
but the argument for that doesn't seem to be immediate.
\end{remark}

$\ $

$\ $

{\bf Acknowledgements.} 

\vspace{4pt}

\noindent
The first-named author would like to thank Marcelo Aguiar and
Frank Sottile for introducing him to combinatorial Hopf algebras.

\vspace{4pt}

\noindent
The second-named author would like to thank
Ed Effros, Sergey Fomin and Dan Voiculescu for inspiring discussions
on the possible use of Hopf algebras in the combinatorics of free
probability. Some of these discussions took place during two workshops
at the American Institute of Mathematics in Palo Alto, in January 2005 
and in June 2006; the second-named author gratefully acknowledges his 
participation in these workshops.

$\ $

$\ $

$\ $

Mitja Mastnak

Department of Pure Mathematics, University of Waterloo.

\vspace{4pt}

Current address: Department of Mathematics and Computer Science,

Saint Mary's University,

Halifax, Nova Scotia B3H 3C3, Canada.

Email: mmastnak@cs.smu.ca

$\ $

Alexandru Nica

Department of Pure Mathematics, University of Waterloo,

Waterloo, Ontario N2L 3G1, Canada.

Email: anica@math.uwaterloo.ca

\end{document}